\documentclass[journal]{IEEEtran}
\IEEEoverridecommandlockouts

\usepackage{amsmath,amssymb}
\usepackage{graphicx}
\usepackage{subcaption}
\usepackage{cite,url}
\usepackage{psfrag}
\usepackage{enumerate}
\usepackage{algorithmic}
\usepackage{algorithm}
\usepackage{acronym}
\usepackage{color, soul}
\usepackage{amsthm}
\usepackage{pifont}

\newcommand{\refeq}[1]{(\ref{#1})}
\newcommand{\reffig}[1]{Figure~\ref{#1}}

\newcommand{\reftab}[1]{Table~\ref{#1}}

\newcommand{\refsec}[1]{Section~\ref{#1}}

\newcommand{\xmark}{\ding{55}}%

\newtheorem{definition}{Definition}

\newtheorem{remark}{Remark}
\begin{document}

\acrodef{VRS}{Vortex Ring State}
\acrodef{WBS}{Windmill-Brake State}
\acrodef{TWS}{Turbulence Wake State}
\acrodef{NED}{North-East-Down}
\acrodef{PWM}{Pulse Width Modulation}
\acrodef{PMP}{Pontryagin Maximum Principle}
\acrodef{NLP}{Nonlinear Programming}

%
\title{Optimal-time Quadcopter Descent Trajectories Avoiding the Vortex Ring and Autorotation States}
\author{Amin~Talaeizadeh$^a$, Duarte Antunes$^{b,*}$, Hossein Nejat Pishkenari$^{a,*}$ and Aria Alasty$^a$
\thanks{a  Mechanical Engineering Department, Sharif University of Technology, Tehran, I.R.Iran

b Control Systems Technology Group, Department of Mechanical Engineering, Eindhoven University of Technology, 5600MB Eindhoven, The Netherlands

Email addresses: amtalaiezadeh@mech.sharif.edu (A. Talaeizadeh), d.antunes@tue.nl (D. Antunes), nejat@sharif.edu (H. Nejat Pishkenari), aalasti@sharif.edu (A. Alasty).
}}
\maketitle
\IEEEpeerreviewmaketitle
\begin{abstract}
It is well-known that helicopters descending fast may enter the so-called \ac{VRS}, a region in the velocity space where the blade’s lift differs significantly from regular regions. This may lead to instability and therefore this region is avoided, typically by increasing the horizontal speed. This paper researches this phenomenon in the context of small-scale quadcopters. The region corresponding to the \ac{VRS} is identified by combining first-principles modeling and wind-tunnel experiments. Moreover, we propose that the so called \ac{WBS} or autorotation region should also be avoided for quadcopters, which is not necessarily the case for helicopters.  A model is proposed for the velocity constraints that the quadcopter must meet to avoid these regions. Then, the problem of designing optimal time descend trajectories that avoid the \ac{VRS} and \ac{WBS} regions is tackled. Finally, the optimal trajectories are implemented on a quadcopter. The flight tests show that by following the designed trajectories, the quadcopter is able to descend considerably faster than purely vertical trajectories that also avoid the \ac{VRS} and \ac{WBS}.
\end{abstract}
\acresetall
\section{Introduction}
\label{sec:Introduction}
\IEEEPARstart{Q}{uadcopters} have been introduced in recent years and are expected to have a large impact on many applications, such as agriculture monitoring, industrial inspection, entertainment industry, among others~\cite{cai2014survey}. One of the interesting features of these aerial robots is that they are able to perform aggressive maneuvers due to their agility. In order to safely conduct these aggressive maneuvers, it is important to model and understand the dynamic behavior at high speeds, where non-trivial aerodynamic effects come into play. 

Such aerodynamic effects have been extensively studied in the context of helicopters, where one can find several seminal papers~\cite{saunders1975dynamics, johnson2005model, peters1982momentum, basset2008prediction, leishman2006principles, seddon2001basic}.  As explained in these papers, there are three regions in the velocity space called \ac{VRS}, \ac{WBS} and \ac{TWS} in which the behavior of the helicopter and the effect of the environment is different from regular regions. The \ac{VRS} and \ac{TWS} regions must be avoided, but it is still possible, although undesired, to control the helicopter in the \ac{WBS} region. However, these effects have received only limited attention in the context of quadcopters, despite the fact that they can play a very important role due to their agility. Some exceptions are~\cite{huang2009aerodynamics, westbrook2014investigation, foster2017high, hoffmann2011precision, chenglong2015vortex, jimenez2001induced, prasad2006prediction, taamallah2010qualitative, bangura2012nonlinear}.

In~\cite{huang2009aerodynamics} it is reported that for descending trajectories and for some velocity state space regions a quadcopter can become unstable. Other studies show that these instabilities are due to the \ac{VRS} effects, see~\cite{westbrook2014investigation, foster2017high, hoffmann2011precision}, where simple models are defined for explaining the boundary of these regions. These simple models are similar to helicopter models and are used in some works to design controllers that avoid entering the \ac{VRS}~\cite{chenglong2015vortex, jimenez2001induced, prasad2006prediction, taamallah2010qualitative, bangura2012nonlinear}. In fact, these studies are similar to helicopter \ac{VRS} avoiding controllers introduced in~\cite{abildgaard2009active, ribera2007helicopter}. However, the results of these models and control strategies have only been shown through simulation and have not been validated through real measurements. Moreover, due to differences in the blades mechanisms between helicopters and quadcopters\footnote{The helicopters’ blades have a constant rotational speed and their pitch angle is variable with a swashplate mechanism for changing the direction of the blade disk thrust. In contrast, quadcopters’ motor speed is variable, and their blades’ pitch angle is usually fixed. Moreover, the hings of the helicopter blades provide a flapping degree of freedom, whereas in contrast to the quadcopter where the blades are rigid.}, high fluctuation regions are different. Furthermore, such controllers are not time-optimal, which is often required in the context of aggressive maneuvers. (see \cite{hehn2012performance, liang2018dynamics, phang2015systems} which tackle the problem of finding time-optimal trajectories without considering aerodynamic constraints.)

In this paper, we conduct static experimental tests on a wind tunnel to identify the regions where fluctuations occur. Combining this information with the \ac{VRS} and \ac{TWS} regions predicted by a theoretical model, we provide a simple model for the region in the velocity space to be avoided by the quadcopter. As we will discuss, we claim that the quadcopter, contrarily to helicopters, should also avoid the \ac{WBS} region. According to the model for constraints and state bounds, different optimal trajectories are designed, considering planar 2D and full 3D models for the motion of the quadcopter. Due to the complexity of this nonlinear optimal control problem, we resort to the optimal control toolbox, GPOPS-II~\cite{patterson2014gpops}. The designed trajectories to descend as fast as possible boil down to intuitive maneuvers, such as oblique and flips maneuvers in 2D, and helix type maneuvers in 3D. Finally, the optimal trajectories are implemented on a quadcopter. The flight tests show that by following the designed trajectories, the quadcopter is able to descend considerably faster than purely vertical trajectories that also avoid the \ac{VRS}.

The remainder of the paper is organized as follows. In \refsec{sec:Problem Formulation}, the aerodynamic effects in descent maneuvers that may cause instability are described, and two major regions in the descent of quadcopters are introduced, based on a theoretical model borrowed from helicopter aerodynamic theory. In \refsec{sec:Identification}, the wind tunnel test results are described and a model is suggested for the boundaries of the unstable regions in the velocity space. \refsec{sec:Optimal Trajectories} provides optimal minimum time trajectories in 2D and 3D spaces that are designed considering the \ac{VRS} constraints, a model of the quadcopter, and GPOPS-II. \refsec{sec:Experiments} presents the flight test experiments by which the quadcopter follows the optimal time trajectories while avoiding the \ac{VRS}. \refsec{sec:Conclusions} provides a brief discussion and research lines for future work. 

\section{Model of Critical Regions During Descent}
\label{sec:Problem Formulation}
\begin{figure}[tbp]
\centering
\includegraphics[trim={0.5cm 0.5cm 0.5cm 0.5cm},clip, width=.9\linewidth]{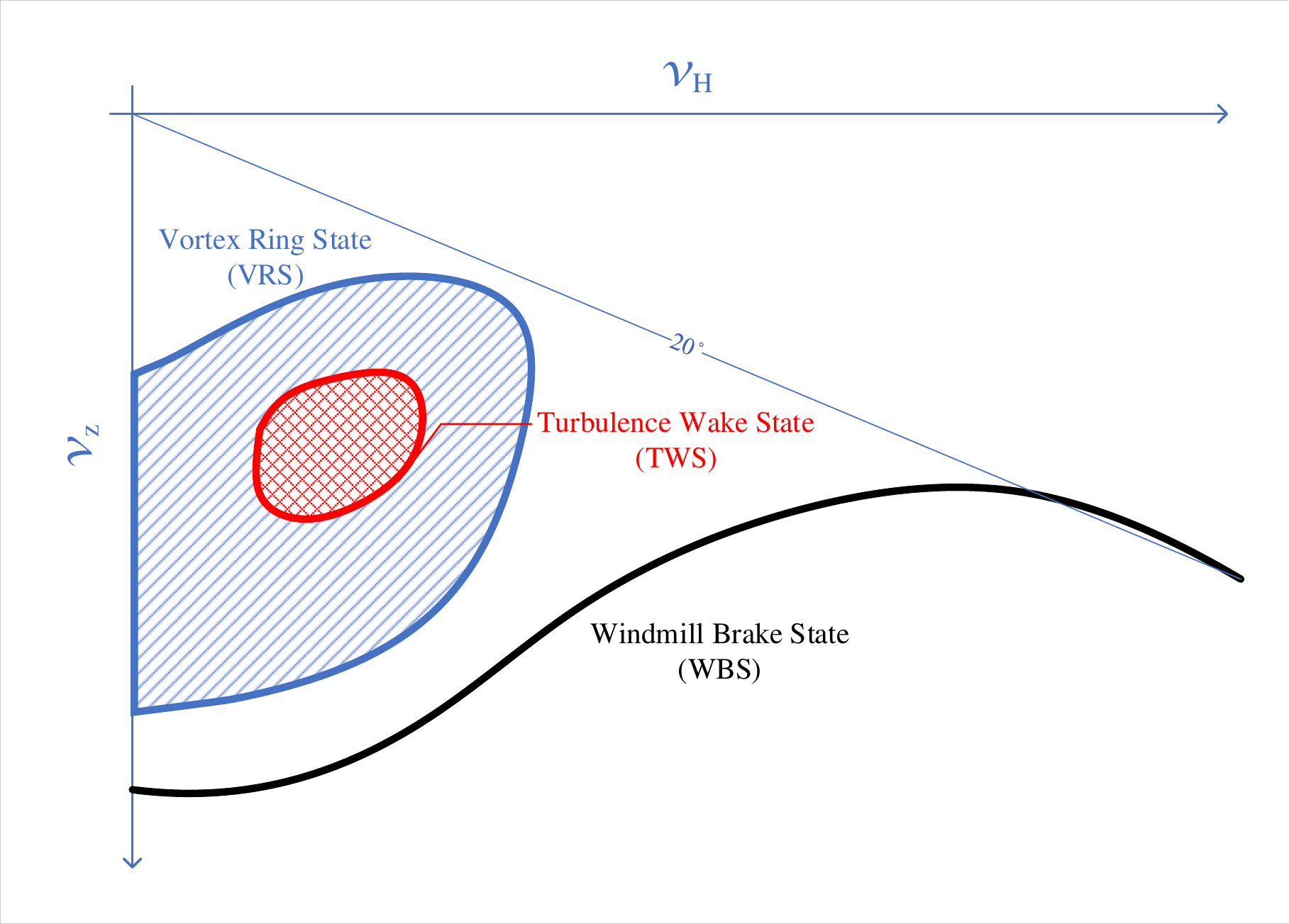}
\caption{The \ac{VRS}, \ac{TWS} and \ac{WBS} regions in the velocity state space for helicopters are shown in the figure. Both vertical, $v_z$, and horizontal $v_H$ speeds are normalized by the induced velocity at hover, where $v_H$ and $v_z$ are horizontal and vertical speeds in the blade disk frame, respectively~\cite{saunders1975dynamics}.}
\label{fig:VRSHelicopters}
\end{figure}

%
In this section, the behavior of the quadcopter in descent maneuvers is described based on the theory available for helicopters. When a helicopter is descending, in some regions in the velocity state space, intensive vibrations and fluctuations in thrust occur. By increasing the horizontal speed of the helicopter, these effects can be reduced. Depending on the vertical and horizontal speeds, the amplitude of fluctuations in thrust can vary. The low and high amplitude fluctuations are known as \ac{VRS} and \ac{TWS}, respectively. In \reffig{fig:VRSHelicopters}, the \ac{VRS} and \ac{TWS} regions are illustrated in the velocity space normalized by the induced velocity at hover. The induced velocity corresponds to the flow at the tip of the blade from the high-pressure air at the bottom of the blade to the low-pressure air at the top of the blade. 
For a helicopter, the descent modes can be divided into four modes:
\begin{itemize}
\item \emph{Normal Descent:} In this mode, there is no fluctuation on the thrust. This corresponds to a normal flight mode.
\item \emph{\ac{VRS}:} In a descent maneuver, the first low fluctuations in the thrust are in the Vortex Ring state region. This region is described in \refsec{subsec:VRS}.
\item \emph{\ac{TWS}:} In this region, because of high amplitude fluctuations on the thrust, the helicopter experiences intensive vibrations which might lead to instability. This region is also considered as \ac{VRS} in this paper because this region is encompassed by the \ac{VRS} region in the velocity state space.
\item \emph{\ac{WBS}:}(or autorotation) In this region, while the helicopter still suffers high amplitude fluctuations, there is enough lift created by the airflow (contrarily to the \ac{VRS}) such that the helicopter can still be controlled.\footnote{A helicopter typically enters \ac{WBS} or autorotation region when it cannot use its engine because of motor's or tail-rotor's malfunctioning.}  This can be achieved by switching off the motors and changing the pitch angle of the blades to simply control the angle of attack (the pilot can, e.g., change a positive angle of attack to a negative one to produce the thrust and change the thrust vector direction with the swashplate to control the orientation of the helicopter~\cite{seddon2001basic}). The blades lift (and rotation) will now be caused by the air moving up through the blade disk and not by the motors. This effect is equivalent to the gliding flight in fixed-wing planes. In \refsec{subsec:WBS}, the differences between helicopter and quadcopter, and the behavior of the quadcopter in this region are explained.

\end{itemize}
To determine the behavior of the quadcopter in the descent phases, firstly, we will provide a model for the \ac{VRS} and \ac{TWS} regions. For the \ac{WBS} region, we will discuss that due to the differences between the helicopter's and quadcopter's blade disks, the \ac{WBS} models for helicopters are not valid for quadcopters, and we will find the behavior of the quadcopter's blade disks by wind tunnel tests which are provided in \refsec{sec:Identification}.

\subsubsection{Vortex Ring State Model in Quadcopter}
\label{subsec:VRS}
\begin{figure}[tbp]
	\centering
	\includegraphics[trim={0.5cm 1cm 0.5cm 0.5cm},clip, width=.9\linewidth]{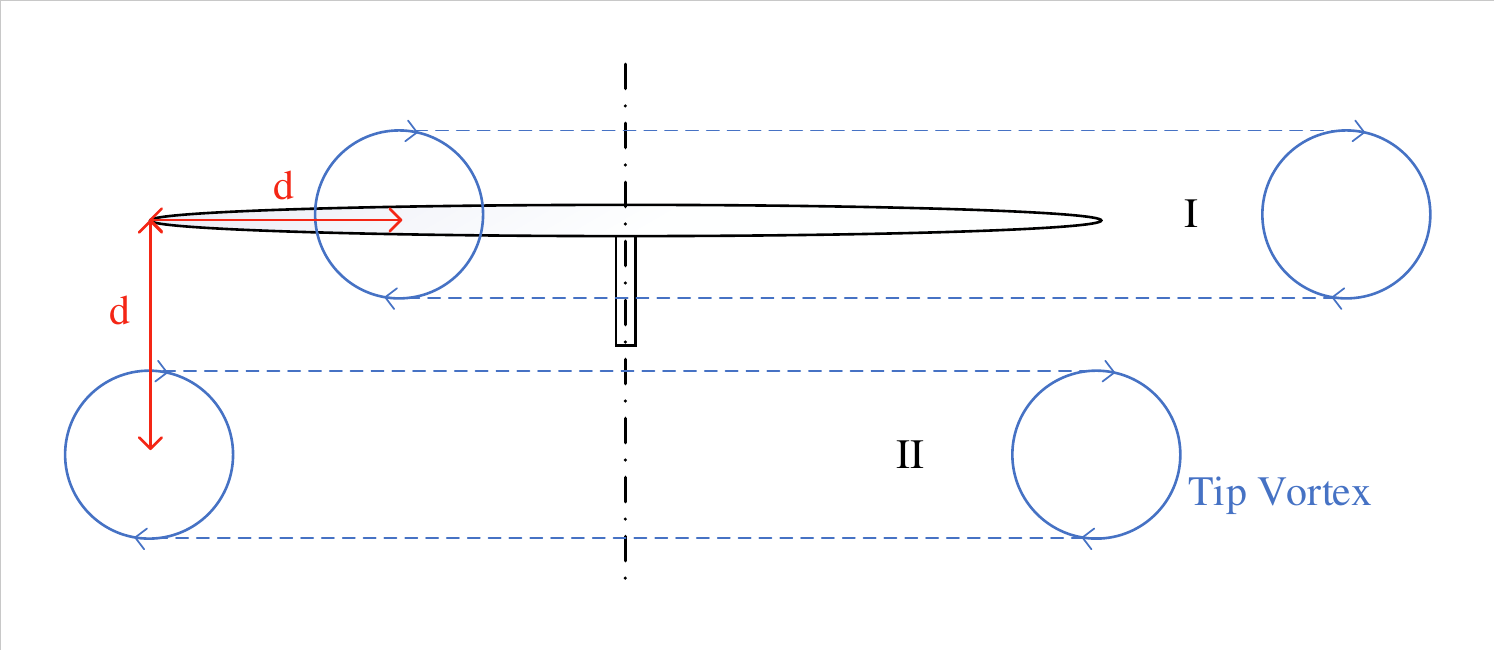}
	\caption{Asymmetry in maneuvers with horizontal speeds pertaining to Onera's criterion for the vortex ring state model~\cite{jimenez2001induced}. I and II are two vortex rings leaving the rotor with the same velocity, but moving in the parallel and perpendicular directions, respectively. After a time $t$, the they will move the distance $d$, however \reffig{fig:Onera2} shows that tip vortex I is still in contact with the blade disk while vortex II is not.}
	\label{fig:Onera2}
\end{figure}

The momentum theory is not valid in the \ac{VRS} and \ac{TWS} regions. However, it is possible to find the boundary of these regions to identify if the quadcopter is in these regions or not. For modeling the \ac{VRS} and \ac{TWS} boundaries, there are many suggestions in the literature~\cite{castles1951empirical, yaggy1963wind, azuma1966experiments}. One of the closest models to the experiments is Onera’s model~\cite{taghizad2002experimental}. Let $\Vert \textbf{\emph{V}}_{tv} \Vert$ be the speed of blade's tip vortex and let $\epsilon_{\text{cr}}$ be the critical velocity of the blade's tip vortexes. Then the following criterion defines the critical region~\cite{taghizad2002experimental}:
\begin{equation}\label{eq:OneraCriteria}
\Vert \textbf{\emph{V}}_{tv} \Vert \leq \epsilon_{\text{cr}},
\end{equation}
where the speed of blade's tip vortex can be calculated from:
\begin{figure}[tbp]
	\centering
	\includegraphics[trim={0.5cm 1cm 0.5cm 0.5cm},clip, width=.9\linewidth]{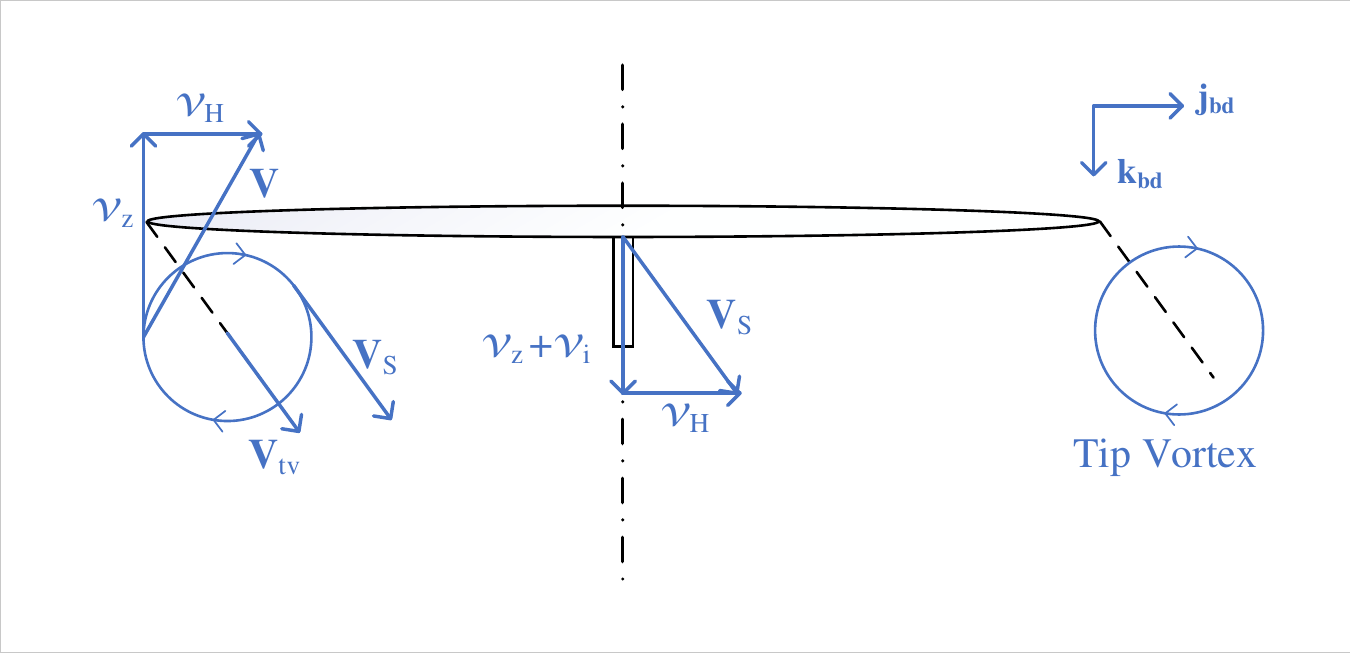}
	\caption{Tip vortexes velocities in oblique descent, used in Onera's criterion~\cite{jimenez2001induced}. $V_S$ in this figure is the velocity of air stream in the blade disk.}
	\label{fig:Onera1}
\end{figure}
\begin{equation}\label{eq:OneraCriteria2}
\textbf{\emph{V}}_{tv}=v_H j_{bd} + (v_i/2+v_z) k_{bd},
\end{equation}
where $v_H$ and $v_z$ are the horizontal and vertical speeds in the blade disk frame, respectively, $v_i$ is the induced velocity in the blade disk, and the unity vectors, $j_{bd}$ and $k_{bd}$ are parallel and perpendicular to the blade disk, as shown in \reffig{fig:Onera1}. Using \refeq{eq:OneraCriteria}, \refeq{eq:OneraCriteria2}, the criterion can be rewritten as:
\begin{equation}
\Vert \textbf{\emph{V}}_{tv} \Vert = \sqrt{ V_{tv_H}^2+V_{tv_z}^2} \leq \epsilon_{\text{cr}}.
\end{equation}
where $V_{tv_H} = v_H$ and $V_{tv_z} = v_i/2+v_z$. \reffig{fig:Onera1} illustrates this setting. However, this model is not symmetric when the quadcopter has a horizontal speed, as shown in \reffig{fig:Onera2}. Let us consider I and II as two vortex rings leaving the rotor with the same velocity, but moving in the parallel and perpendicular directions, respectively. After a time $t$, the they will move the distance $d$, however \reffig{fig:Onera2} shows that tip vortex I is still in contact with the blade disk while vortex II is not.

\begin{figure}[tbp]
	\centering
	\includegraphics[trim={3.5cm 8cm 4.5cm 9cm},clip, width=.9\linewidth]{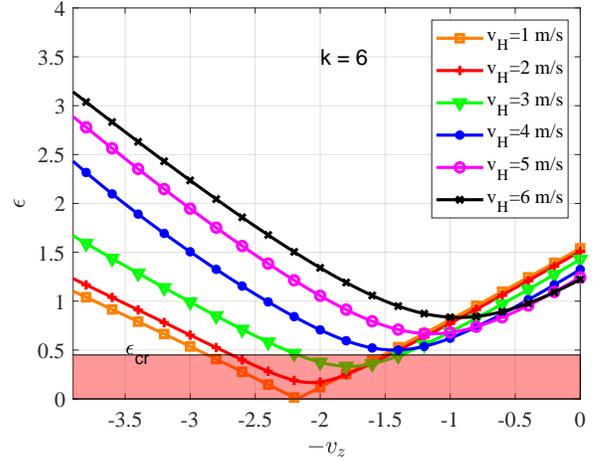}
	\caption{Variation of $\epsilon$ based on different values of $v_H$ and $v_z$. The effect of the horizontal speed in the blade disk frame is obvious, more horizontal speed, more $\epsilon$ value. In this figure, $T_{Hover}=0.5~kg$, $\rho =1.2 kg/m^3$, and the blade disk diameter is $0.2 m$.}
	\label{fig:Onera4}
\end{figure}

To improve the accuracy of this model, assuming asymmetric effects, a correction coefficient, $\text{k}$, should be added to the criterion as introduced in the following equation:
\begin{equation}\label{eq:OneraEq}
\sqrt{(v_H/{\text{k}})^2+(v_i/2+v_z)^2} \leq \epsilon_{\text{cr}},
\end{equation}
where $\text{k}$ is a parameter larger than one, $\text{k} \geq 1$. Then $V_{tv_H} = v_H/{\text{k}}$ and $V_{tv_z} = v_i/2+v_z$. Therefore, in order to meet~\refeq{eq:OneraEq} the following conditions are necessary:
\begin{equation}\label{eq:OneraEqk1}
\left\{\begin{array}{ll}
V_{tv_H} \leq \epsilon_H\\
V_{tv_z} \leq \epsilon_z       
\end{array} \right.
,\epsilon_z \leq \epsilon_H,
\end{equation}
where $\epsilon_H = \text{k} \epsilon_z$ and $ \epsilon_z =  \epsilon_{\text{cr}}$. $\epsilon_z$ and $\epsilon_H$ are critical blade's tip vortex velocity in pure vertical and horizontal maneuvers, respectively. To calculate the induced velocity, the following equation can be used~\cite{peters1982momentum}:
\begin{equation}\label{eq:InducedVel}
v_i \sqrt{v_H^2+(v_i+v_z )^2}=v_h^2,
\end{equation}
where $v_h$ is given constant coinciding with the induced velocity in the blade disk at hover and can be calculated as:
\begin{equation}\label{eq:InducedVelHover}
v_h=\sqrt{T_{\text{Hover}}/2\rho \text{A}},
\end{equation}
where $T_{\text{Hover}}$ is the thrust of the motor at hover, $\rho$ is the air density, and $\text{A}$ is the blade disk area. Usually, $T_{\text{Hover}}/ \text{A}$ is denoted by the disk load. Let
\begin{equation}\label{eq:epsilon}
\epsilon=\sqrt{(v_H/{\text{k}})^2+(v_i/2+v_z)^2}.
\end{equation}
If we fix  $v_H$, we can compute $v_i$ as a function of $v_z$ from \refeq{eq:InducedVel} and plot $\epsilon$, as a function of $v_z$. This plot is shown in \reffig{fig:Onera4} for different values of $v_H$. This figure illustrates that by increasing the horizontal speed in the blade disk frame $v_H$, the $\epsilon$ value will increase and therefore it is possible to increase the vertical (descent) speed $v_z$ without entering the \ac{VRS} and \ac{TWS} regions. Such regions are defined by $\epsilon\leq  \epsilon_{cr}$, and are also depicted in \reffig{fig:Onera4}.
%

\begin{figure}[tbp]
	\centering
	\includegraphics[trim={3.5cm 8cm 4cm 8cm},clip, width=.9\linewidth]{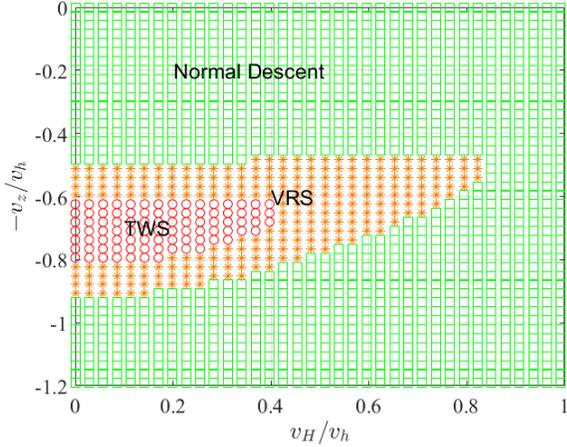}
	\caption{Numerical solution of the Onera’s criterion. The \ac{VRS} and \ac{TWS} are illustrated by orange stars and red circles, respectively. The region identified with green squares corresponds to normal descent in Onera's criterion, however this region may include \ac{WBS} in reality. In this figure, $T_{Hover}=0.5~kg$, $\rho =1.2 kg/m^3$, and the blade disk diameter is $0.2 m$.}
	\label{fig:Onera5}
\end{figure}

\reffig{fig:Onera5} shows the numerical solution of criterion \refeq{eq:OneraEq} with $\text{k}=6$ and $\epsilon_{\text{cr}} =0.4$ for the \ac{VRS} (which illustrates with orange stars) and $\epsilon_{\text{cr}} =0.2$ for the \ac{TWS} (which illustrates with red circles). The induced velocity is calculated from \refeq{eq:InducedVel}. \reffig{fig:Onera5} illustrates the \ac{VRS} and \ac{TWS} based on Onera's criterion. Note that following Onera’s criterion one cannot identify the \ac{WBS} which is a limitation of this theoretical model. Such a region will be identified experimentally in the next section. Since the \ac{TWS} is encompassed by the \ac{VRS}, we will simply consider the union of these regions as a prohibited region and denoted this region simply by the \ac{VRS} region.

\subsubsection{Windmill Brake State in quadcopters}
\label{subsec:WBS}
 As explained before, the control of a helicopter in this region is achieved by changing the pitch. This is possible due to the pitch angle degree of freedom in the helicopters’ blade mechanism. However, for quadcopters, there is no degree of freedom in pitch angles, and there is not a swashplate mechanism to handle the thrust vector in the absence of the motor. Since if we turn off the motors of the quadcopter, we will lose the controllability of the quadcopter, the autorotation effect acts as a brake in the blade disk. In the autorotation mode, by the action of upward air, the blades tend to rotate in the opposite direction, the motor performs torque to turn the blades and this may cause fluctuations on the blade disk. Hence, \ac{WBS} regions should be considered as a forbidden area for quadcopters like the \ac{VRS} and \ac{TWS} regions.
%

\section{Identification of the descend regimes in the quadcopters}
\label{sec:Identification}
\begin{figure}[tbp]
	\centering
	\includegraphics[width=.4\linewidth]{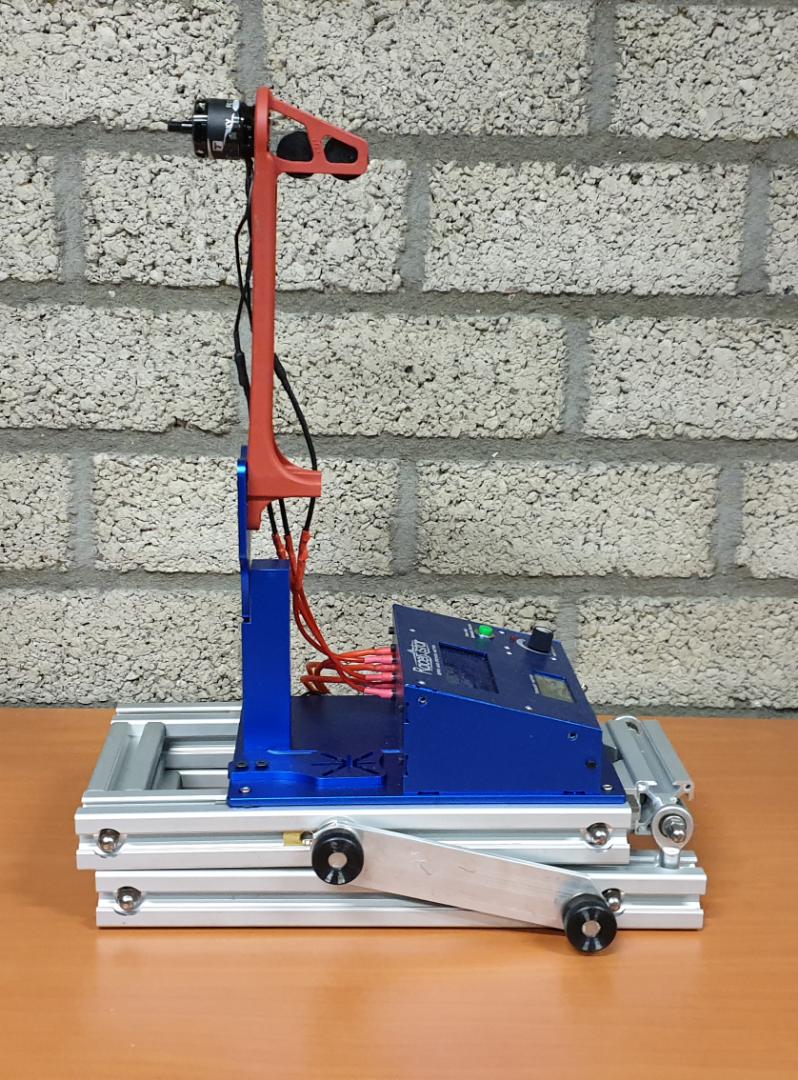}
	\includegraphics[width=.385\linewidth]{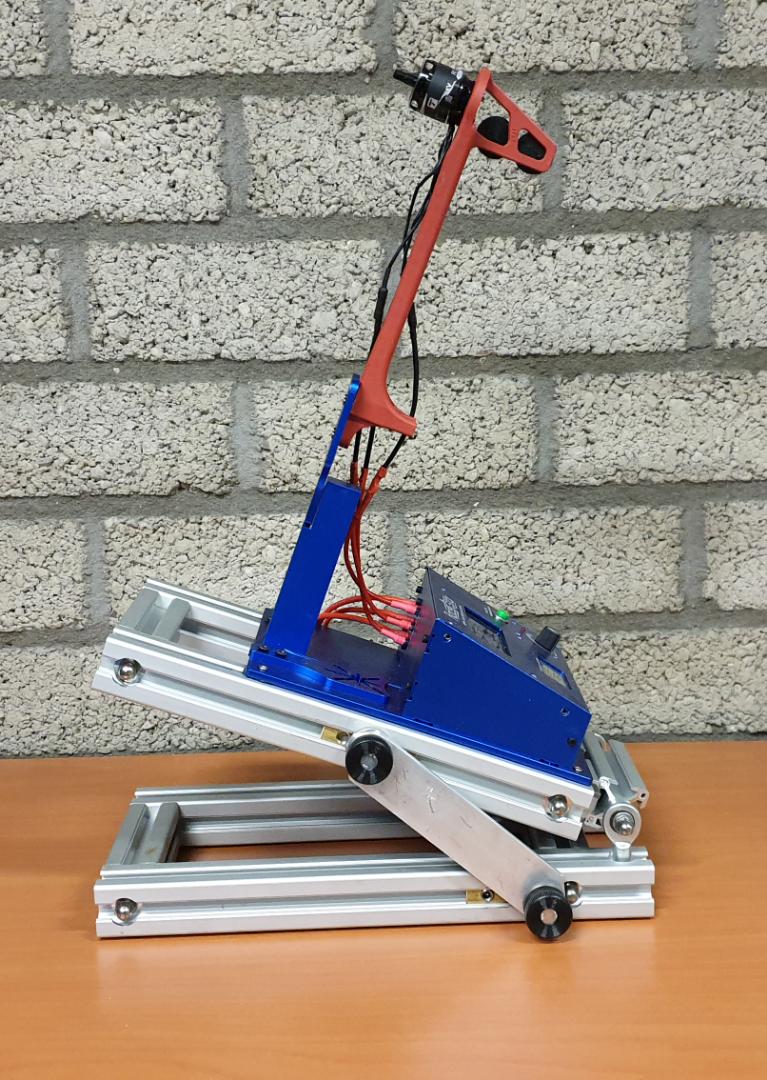}
	\caption{Motor-blade tester to measure thrust, rpm, voltage and \ac{PWM} adapted for simulating various oblique descent angles in wind tunnel. Setup angle is measured by a digital inclinometer.}
	\label{fig:Tester}
\end{figure}

\begin{figure}[tbp]
	\centering
	\includegraphics[trim={0.5cm 0.4cm 0.5cm 0.4cm},clip, width=.9\linewidth]{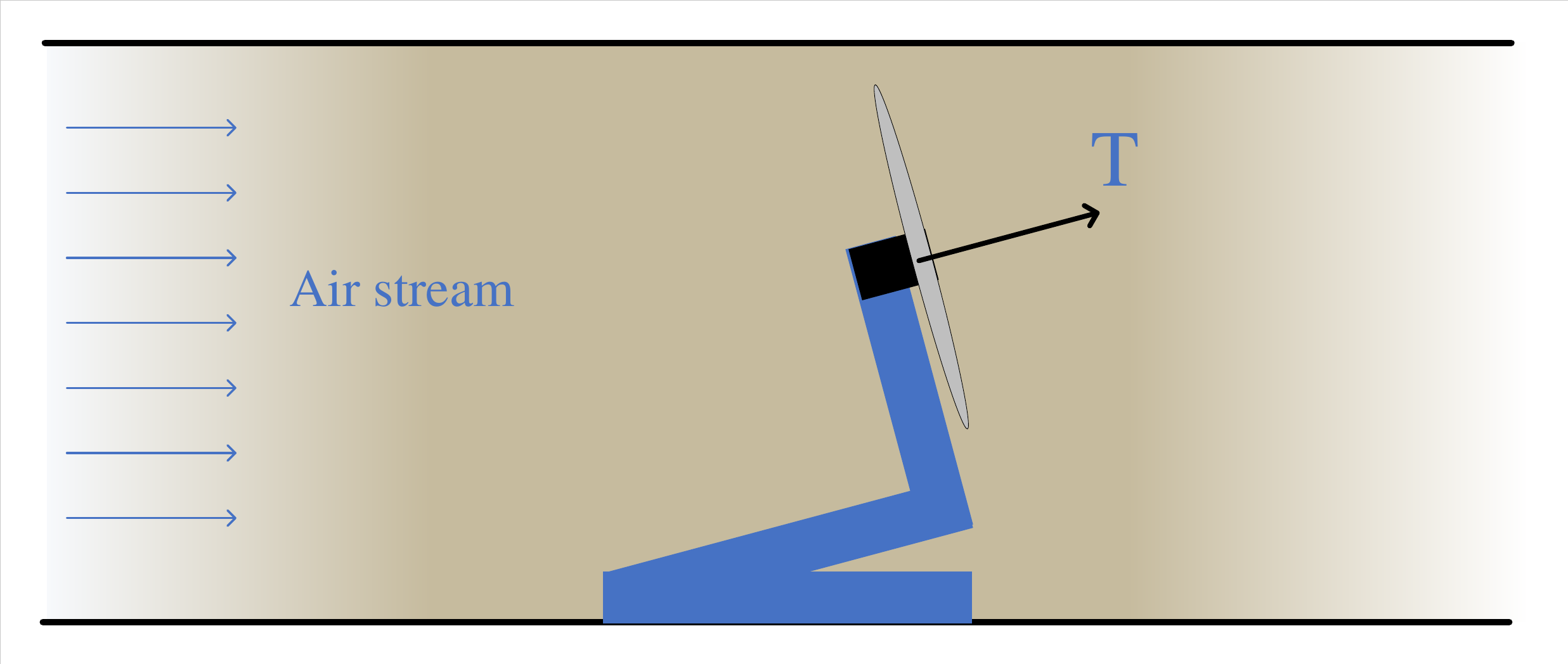}
	\caption{Emulation of descent in wind tunnel with variable angle mechanism for thrust tester setup.}
	\label{fig:TesterWindtunnel}
\end{figure}

\begin{figure}[tbp]
	\centering
	\includegraphics[trim={3.5cm 8cm 4cm 9cm},clip, width=.9\linewidth]{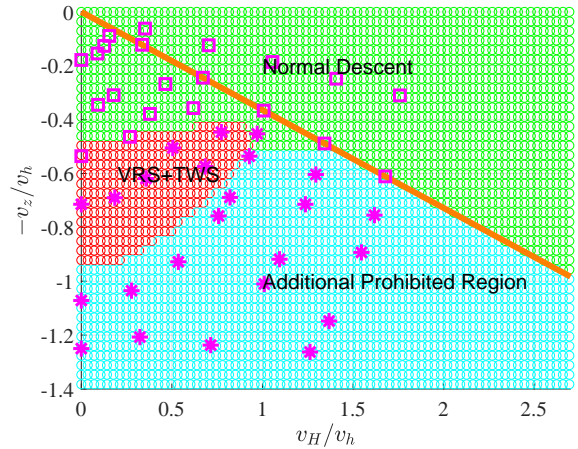}
	\caption{Wind tunnel experiment results where stars and squares define the high fluctuation and low fluctuation points, respectively. The \ac{VRS} region, using the Onera's criterion is illustrated with red color. The blue color region is also prohibited due to the high fluctuations, and it concluded in this research that these fluctuations are the \ac{WBS} effect. Normal regions are depicted by green color. To simplify this constraint, the top region of the orange line is considered as the allowable region}
	\label{fig:VRS_Quad}
\end{figure}

To identify the behavior of the brushless motors with fixed pitch angle blades, in different descent velocities and angles, an experimental setup has been designed, see \reffig{fig:Tester}. This experimental setup can measure the force of motor-blades, rpm, current, voltage. The angle of the setup can be modified so that it emulates the oblique descent in the wind tunnel, see \reffig{fig:TesterWindtunnel}. For these experiments, the wind tunnel of the Eindhoven University of Technology (TU/e) is used. In this wind tunnel, we measured the amount of thrust and its fluctuation, the motor rpm, and the motor \ac{PWM} and voltage input in different wind velocities and different setup angles (to simulate the different angles of descent). These tests were performed with different blades and motors. However, the results, illustrated in~\reffig{fig:VRS_Quad} are for T-Motor Air 2213 brushless motor and T9545 Propeller. The results are normalized by the induced velocity at hover. Regions with high fluctuations ($\Delta T /T \geq 10\%$ for thrust fluctuations) are defined as prohibited regions.

\begin{figure}[tbp]
	\centering
	\includegraphics[trim={3cm 0.5cm 3cm 0.5cm},clip, width=.9\linewidth]{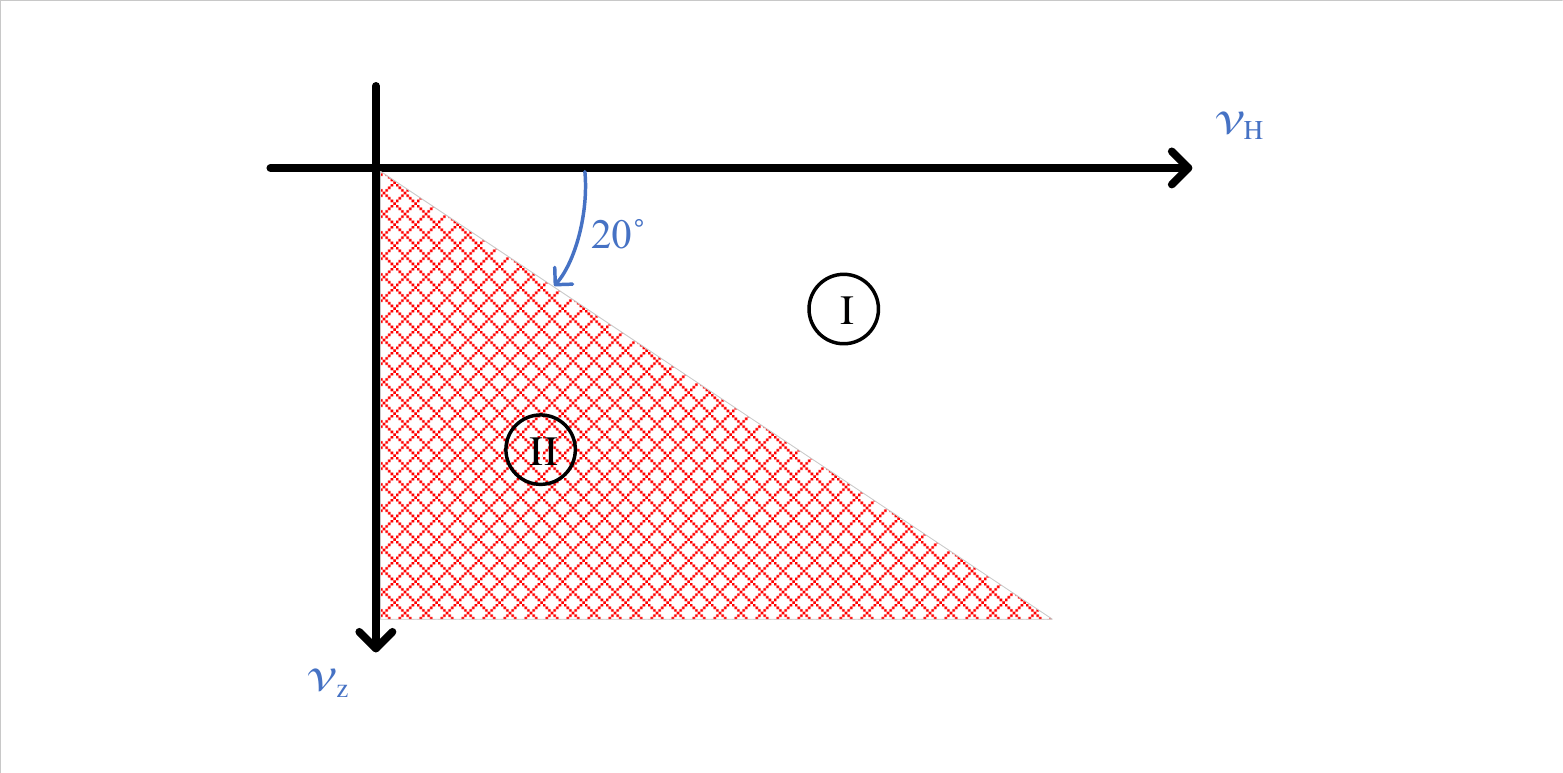}
	\caption{The simplified constraint used for optimal trajectory design. Regions I and II are allowed and prohibited regions, respectively.}
	\label{fig:SimpleConstraint}
\end{figure}

As predicted, based on the theoretical analysis, high fluctuations occur in the \ac{VRS} region computed from the Onera’s criterion. Interestingly, fluctuations also occur in the blue regions shown in \reffig{fig:VRS_Quad}, which is identified as the \ac{WBS} (which as mentioned before, is not predicted from Onera’s model). As also depicted in \reffig{fig:VRS_Quad}, and also in \reffig{fig:SimpleConstraint}, we can define a simplified region where high fluctuations do not occur. This is summarized in the next definition.

%
%
\begin{definition}{Simplified \ac{VRS} Velocity Constraint}
	\begin{equation}\label{eq:Constraint}
	\begin{array}{l}
		|v_{H}| \tan{20^{\circ}} \geq v_{z},\\
		\forall v_{z}\leq 0 , v_H \in \rm I\!R,
	\end{array}
	\end{equation}
where $v_H$ and $v_z$ are in the body fixed frame.
\end{definition}
\begin{remark}
	The motion planner should design a trajectory to avoid these regions. Note that to take into account these constraints it is necessary to compute the relative wind velocities in the blade disk frame, not in the inertial frame, as it will be discussed in the next section.
\end{remark}

\begin{remark}
	The dimensions of the velocity constraints are normalized by the induced velocity at hover. Therefore, this model can be used for every kind of quadcopters' blade disks with different disk loads and diameters that have fixed-pitch blades.
\end{remark}

\section{Design of Optimal Trajectories}
\label{sec:Optimal Trajectories}
This section tackles the problem of designing optimal trajectories for a quadcopter to descend as fast as possible, considering the \ac{VRS} and \ac{WBS} constraints. \refsec{subsec:QuadEoM} provides the equations of motion of a quadcopter. In~\refsec{sec:2Danal} and~\refsec{sec:3Danal}, the motion of the quadcopter is considered in 2D and 3D space, respectively, and descend trajectories with minimum time are designed.

\subsection{Quadcopter Equations of Motion}
\label{subsec:QuadEoM}
If we consider the \ac{NED} as the inertial frame, and a body fixed frame, denoted by $\text{b}$ with its origin coinciding with the center, z axis aligned with the gravitational vector when the quadcopter is at hover and $x$ and $y$ axis as shown in \reffig{fig:QuadFrame}, the total external force $\textit{\textbf{F}}^{\text{I}}$ expressed in the inertial frame acting on the quadcopter can be expressed as follows (see \cite{talaeizadehaccurate}):
\begin{figure}[tbp]
	\centering
	\includegraphics[trim={3cm 0.5cm 1cm 0.5cm},clip, ,width=.9\linewidth]{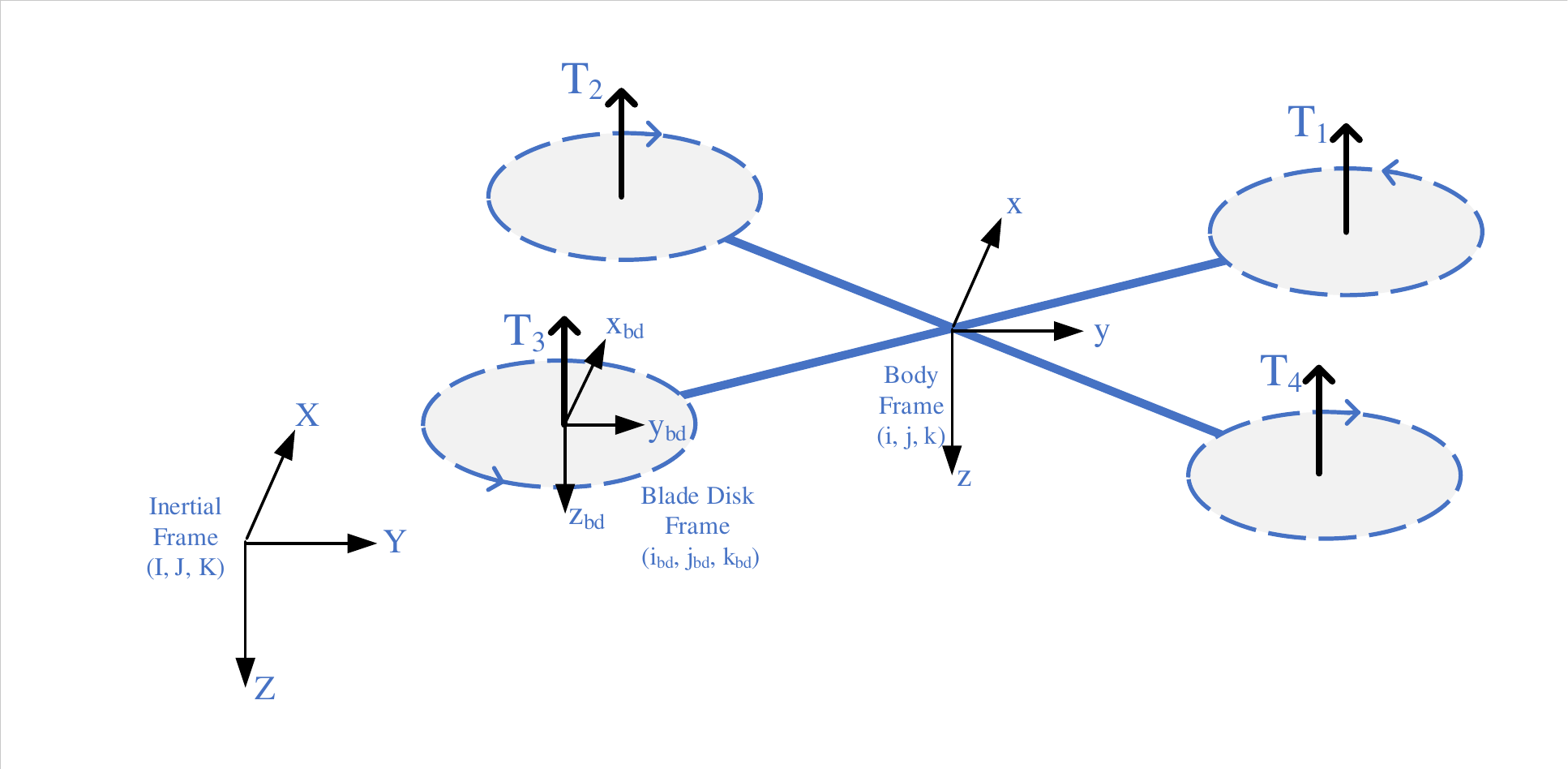}
	\caption{Quadcopter, its motors' thrusts, the inertial, body and blade disk \ac{NED} frames are illustrated.} 
	\label{fig:QuadFrame}
\end{figure}
\begin{equation} \label{eq:DOA}
\sum \textbf{\textit{F}}^\text{I}=\textbf{\emph{R}}_{\text{b}}^{\text{I}}\left\{\begin{array}{c}
0\\
0\\
\sum_{1}^{4}(-T_i)
\end{array} \right\}
+ \left\{\begin{array}{c}
0\\
0\\
\text{mg}
\end{array} \right\}
+\textbf{\textit{F}}_{\text{ext}},
\end{equation}
where $\textbf{\emph{F}}_{ext}$ is disturbance force vector, $\textbf{\emph{R}}_{\text{b}}^{\text{I}}$ is the rotation matrix from the body to the inertial frame and $T_i$ denotes the thrust generated by the $i\textsuperscript{th}$ motor.

By the Euler-Newton formalism, the summation of torques that are applied to the quadcopter is equal to the time rate of change of the angular momentum.
\begin{equation}\label{eq:Newton2}
\sum \textbf{\textit{M}}=d\textbf{\textit{H}}/dt=\dot{\textbf{\textit{H}}}+ \boldsymbol{\omega} \times \textbf{\textit{H}}.
\end{equation}
The torques that are applied to the robot from motors depend on the configuration of the quadcopter. For an X shaped quadcopter, the torques are given by: 

\begin{equation}\label{eq:XshapeTorques}
\textbf{\textit{M}}=
\left\{\begin{array}{c}
M_x\\
M_y\\
M_z
\end{array} \right\}
=
\left\{\begin{array}{c}
b(-T_1+T_2+T_3-T_4)\\
b(T_1+T_2-T_3-T_4)\\
(-M_1+M_2-M_3+M_4)
\end{array} \right\}
+\textbf{\textit{M}}_{\text{ext}},
\end{equation}
where $T_i$ is the output thrust of the $i^{th}$ motor, $b$ is the distance between motors and center of gravity, $M_i$ is the output torque of $i^{th}$ motor and $\textbf{\emph{M}}_{\text{ext}}$ is the external torque vector. In~\refeq{eq:Newton2}, $\emph{\textbf{H}}$ is the angular momentum of the robot that can be calculated as:
\begin{equation}\label{eq:AngularMomentum}
\textbf{\textit{H}}=\textbf{\emph{I}}\boldsymbol{\omega},
\end{equation}
where $\textbf{\emph{I}}$ is the body’s inertia tensor of the quadcopter, and $\boldsymbol{\omega}$ is the angular velocity of the quadcopter’s body.
\begin{equation}\label{eq:omega}
\boldsymbol{\omega}=
\left\{\begin{array}{c}
\omega_x\\
\omega_y\\
\omega_z
\end{array} \right\}
=
\left\{\begin{array}{c}
p\\
q\\
r
\end{array} \right\}.
\end{equation}
Therefore, equation~\refeq{eq:Newton2} can be rewritten as:
\begin{equation}\label{eq:Newton3}
\sum \textbf{\textit{M}}=\textbf{\emph{I}}\dot{\boldsymbol{\omega}}+\textbf{\emph{S}}(\boldsymbol{\omega})\textbf{\emph{I}}\boldsymbol{\omega},
\end{equation}
where
\begin{equation}\label{eq:scewsym}
\textbf{\emph{S}}(\boldsymbol{\omega})=
\begin{bmatrix}
0 & -r & q\\
r & 0 & -p\\
-q & p & 0
\end{bmatrix}.
\end{equation}
Finally, the equations of motion are given by:
\begin{equation}\label{eq:eom}
\begin{array}{l}
\dot{\textbf{\emph{P}}} = \textbf{\emph{V}},\\
\dot{\textbf{\emph{V}}} = \sum \textbf{\textit{F}}^\text{I}/\text{m},\\
\dot{\textbf{\emph{R}}}^{\text{I}}_{\text{b}}= {\textbf{\emph{R}}}^{\text{I}}_{\text{b}}S(\omega),\\
\textbf{\emph{I}}\dot{\boldsymbol{\omega}} = \sum \textbf{\textit{M}}-\textbf{\emph{S}}(\boldsymbol{\omega})\textbf{\emph{I}}\boldsymbol{\omega},\\
\end{array}
\end{equation}
where $\textbf{\emph{P}}$ and $\textbf{\emph{V}}$ are the position and the velocity of the quadcopter's center of mass in the inertial frame, respectively, and $\text{m}$ is the quadcopter's mass.

Alternatively to the equation $\dot{\textbf{\emph{R}}}^{\text{I}}_{\text{b}}= {\textbf{\emph{R}}}^{\text{I}}_{\text{b}}S(\omega)$, we can write:
\begin{equation}
\dot{q}=\textbf{\emph{G}}(q)\boldsymbol{\omega},
\end{equation}
where $q = [q_0, q_1, q_2, q_3]$ are the quaternions and $\textbf{\emph{G}}(q)$ is a matrix of quaternions \cite{baruh1999analytical}.

\subsection{2D Optimal Trajectories}\label{sec:2Danal}
Consider the quadcopter is in the Y-Z plane, as shown in \reffig{fig:2DQuad}, its motion is planar and there are only three degrees of freedom. Let us combine two pairs of motors on each side, so that we can assume a single thrust force is applied on each side. Although this leads to a simple planar model, it is still an underactuated system, such as the full 3D quadcopter model. 
The equations of motion for this simplified can be derived from the full model~\refeq{eq:eom} and are given by:
\begin{figure}[tbp]
	\centering
	\includegraphics[trim={0.5cm 0.5cm 1cm 0.5cm},clip, width=.9\linewidth]{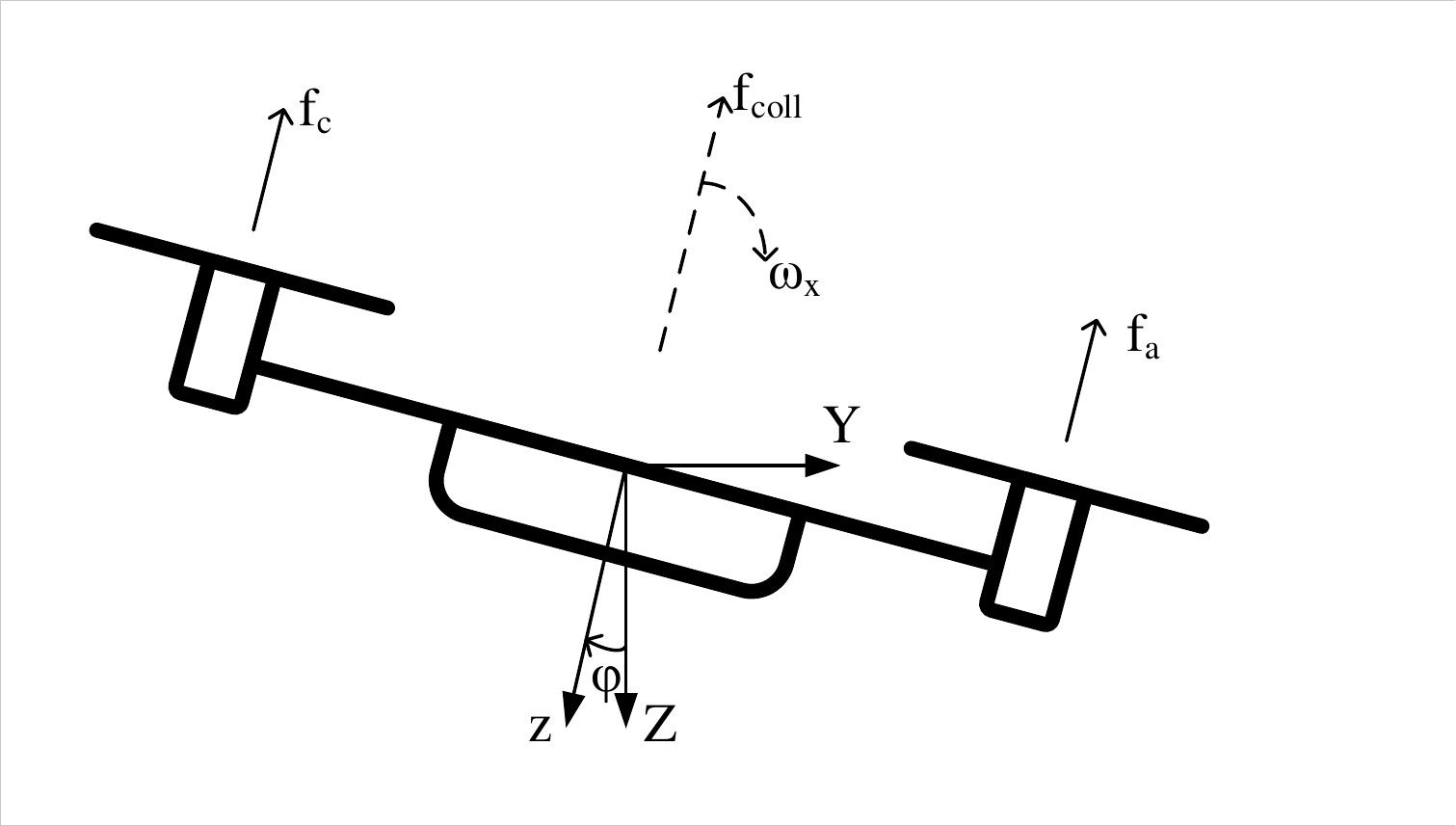}
	\caption{Quadcopter in 2D analysis, including two motors. The inertial and body frames are illustrated.}
	\label{fig:2DQuad}
\end{figure}
\begin{equation}\label{eq:2DEoM2}
\left. \begin{array}{ll}
\ddot{Z} =-f_{\text{coll}} \cos{\phi} + \text{g},\\
\ddot{Y} = f_{\text{coll}} \sin{\phi},\\
I_{xx} \ddot{\phi} = ml(f_a-f_c),
\end{array} \right.
\end{equation}
where $f_{coll}$ is the collective acceleration from the motors' thrusts:
\begin{equation}
f_{coll} = f_a+f_c,
\end{equation}
where $f_a$ and $f_c$ are the thrusts of each motor, divided by the mass. If we assume the angular velocity as a control input, the kinematic equation can be assumed to be simply given by:
\begin{equation}\label{eq:2DEoM}
\left. \begin{array}{ll}
\dot{\phi} = \omega_x.
\end{array} \right.
\end{equation}
The system states, $S$, and inputs, $U$, are:
\begin{equation}
\left. \begin{array}{ll}
S =(Y,V_Y,Z,V_Z,\phi),\\
U = (f_{coll}, \omega_x).\\
\end{array} \right.
\end{equation}
The state space equation is:
\begin{equation}\label{2Deom}
\begin{array}{r}
\begin{bmatrix}
\dot{Y}\\
\dot{V}_Y\\
\dot{Z}\\
\dot{V}_Z\\
\dot{\phi}
\end{bmatrix}
=
\begin{bmatrix}
0 & 1 & 0 & 0 & 0 \\
0 & 0 & 0 & 0 & 0 \\
0 & 0 & 0 & 1 & 0 \\
0 & 0 & 0 & 0 & 0 \\
0 & 0 & 0 & 0 & 0 \\
\end{bmatrix}
\begin{bmatrix}
Y\\
V_Y\\
Z\\
V_Z\\
\phi
\end{bmatrix}
+                            \\
\begin{bmatrix}
0 & 0\\
\sin{\phi} & 0\\
0 & 0\\
-\cos{\phi} & 0\\
0 & 1
\end{bmatrix}
\begin{bmatrix}
f_{\text{coll}}\\
\omega_x
\end{bmatrix}
+
\begin{bmatrix}
0\\ 0\\ 0\\ \text{g}\\ 0
\end{bmatrix}
\end{array}
\end{equation}
Our goal is to design a control law for the quadcopter to move from the initial state at height $h_1$ to the final state at height $h_2$ both at hover.

The cost function for the optimal trajectory design, $J$, is for the minimum time problem:
\begin{equation}\label{eq:cost}
J = \min \int_{0}^{t_f} g(S,U) dt = \min \int_{0}^{t_f} 1 dt.
\end{equation}
Taking into account the relation between velocities in the body fixed frame $(v_H, v_z)$ and the inertial frame $(V_Y, V_Z)$:
\begin{equation}
\left(
\begin{matrix}
v_H \\
v_z
\end{matrix}
\right)
=
\begin{bmatrix}
\cos{\phi} & -\sin{\phi} \\
\sin{\phi} & \cos{\phi}
\end{bmatrix}
\left(
\begin{matrix}
V_Y \\
V_Z
\end{matrix}
\right),
\end{equation}
we can consider the simplified velocity constraint \refeq{eq:Constraint}, $C(V_Y,V_Z,\phi) \geq 0$, where:
\begin{equation}\label{eq:constraintHamilton}
C(V_Y,V_Z,\phi) = V_Z (1- \text{B} \tan{\phi})-V_Y (\text{B}+ \tan{\phi}),
\end{equation}
and $\text{B} = \tan(20^\circ)$. This constraint is obtained by finding the air relative velocity in the body fixed frame and using the Definition~1.
To take into account this constraint we modify the Hamiltonian to:
\begin{equation}
\left. \begin{array}{lll}
H(S,U,\lambda)&=&g(S,U)+\lambda^T f(S,U)+\lambda_C C(V_Y,V_Z,\phi),\\
\end{array} \right.
\end{equation}
%
where $\lambda$ denotes the co-state variables:
\begin{equation}
\lambda = \begin{bmatrix}
\lambda_Y & \lambda_{V_Y} & \lambda_Z & \lambda_{V_Z} & \lambda_{\phi}
\end{bmatrix}^T.
\end{equation}
Using cost function, dynamic and path equations, the Hamiltonian can be written as: 
\begin{equation}\label{eq:Hamilton}
\left. \begin{array}{lll}
H(S,U,\lambda)&=&1+\lambda_Y V_Y+\lambda_{V_Y}  (f_{\text{coll}}  \sin{\phi} )+\lambda_Z V_Z\\
& & +\lambda_{V_Z}  (-f_{\text{coll}}  \cos{\phi}+g)+\lambda_{\phi} \omega_x+\\
& & \lambda_C (V_Z (1-\text{B} \tan{\phi})-V_Y (\text{B}+ \tan{\phi})).\\
\end{array} \right.
\end{equation}
Using the Maximum Principle for problems with mixed inequality constraints conditions~\cite{sethi2019optimal} and Hamitonian in~\refeq{eq:Hamilton}, we obtain:
\begin{equation}
\left. \begin{array}{lll}
\dot{\lambda}(t)&=&-\frac{dH}{dS}(S^* (t),U^* (t),\lambda(t))^T,\\
\dot{\lambda}_Y&=&0 \rightarrow \lambda_Y=c_1,\\
\dot{\lambda}_{v_Y}&=&-\lambda_Y+\lambda_C (\text{B}+\tan{\phi}), \\ 
\dot{\lambda}_Z&=&0 \rightarrow \lambda_Z=c_3,\\
\dot{\lambda}_{v_Z}&=&-\lambda_Z-\lambda_C (1- \text{B} \tan{\phi}),  \\
\dot{\lambda}_{\phi}&=&f_{\text{coll}} (-\lambda_{v_Y}  \cos{\phi}-\lambda_{v_Z} \sin{\phi})\\
& & -\lambda_C/(\cos{\phi})^2 (V_Z \text{B}-V_Y),
\end{array} \right.
\end{equation}
where $\lambda_C \geq 0$ and $\lambda_C C(V_Y,V_Z,\phi) = 0$ are complementary slackness conditions, and $c_1$ and $c_3$ are constants.

This problem is very complex and it is hard to solve this nonlinear optimal control problem analytically, so we use a numerical method to solve it. Namely, we used GPOPS-II, a Matlab plugin, to solve this optimization problem. GPOPS-II implements a new class of variable-order Gaussian quadrature methods where a continuous-time optimal control problem is approximated as a sparse \ac{NLP} problem~\cite{patterson2014gpops}.


Besides the constraints corresponding to the dynamics~\refeq{2Deom}, the constraint that should be satisfied during the whole path is~\refeq{eq:constraintHamilton}.

\begin{figure}[tbp]
	\centering
	\includegraphics[trim={4.5cm 8cm 4.5cm 8cm},clip ,width=\linewidth]{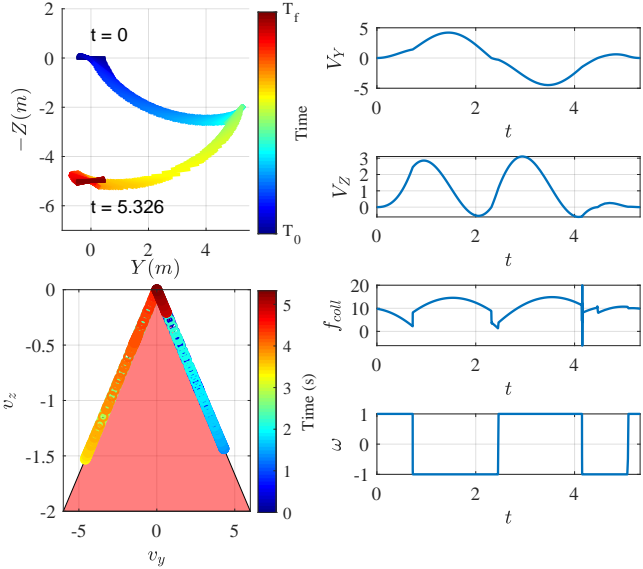}
	\caption{Optimal 2D trajectory design for equal initial and final $Y$ position is shown for the whole trajectory with blue to red colors. The velocity constraint, considering the velocities in the blade disk is satisfied.}
	\label{fig:YfixAllinOne}
\end{figure}

Initial, final and admissible intervals of the system variables are defined in \reftab{table:ICBound}. According to these initial and final values and bounds, the minimum time trajectory can be designed. With this trajectory, the quadcopter can descend 5 meters in 5.33 seconds. The simulation results and designed optimal trajectory are shown in \reffig{fig:YfixAllinOne}. This figure illustrates that even if the final horizontal displacement is not desired, the quadcopter should move in the horizontal direction meanwhile and come back to its initial $Y$, in order to satisfy the velocity constraint. It shows that the pure descent (without any horizontal displacement during the trajectory) is not the fastest trajectory to descend. 

\begin{figure}[tbp]
	\centering
	\includegraphics[trim={4.5cm 8cm 3.5cm 7cm},clip ,width=\linewidth]{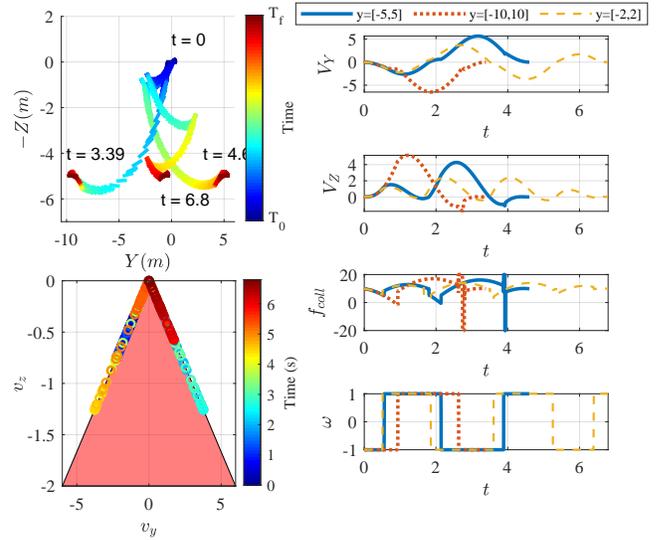}
	\caption{Optimal 2D trajectories design for free (final $Y$ is not fixed but it is bounded) and bounded horizontal displacement, that let the system have oblique and Zig-Zag trajectories, respectively. Moreover, the $\phi$ interval is $[-\pi/3, \pi/3]$. It is obvious that if the pure oblique descent is allowed, then the trajectory will be faster than Zig-Zag maneuver. Moreover, by increasing the horizontal displacement limitation, the trajectory duration will be increased.}
	\label{fig:YfreeAllinOne}
\end{figure}

\begin{table}[tbp]
\centering
\caption{Initial and final values as well as admissible intervals of the system variables}
\begin{tabular}{cccc}
\hline \noalign{\smallskip} 
\textbf{Variable} & \textbf{Initial Cond.} & \textbf{Final Cond.} & \textbf{Admissible Intervals} \\\hline\noalign{\smallskip}
$Y$ & 0 & 0 & $[-15, 15]$ \\
$V_Y$ & 0 & 0 & $[-10, 10]$ \\
$Z$ & 0 & $5$ & $[-15, 15]$ \\
$V_Z$ & 0 & 0 & $[-10, 10]$ \\
$\phi$ & 0 & 0 & $[-\pi/3, \pi/3]$ \\
$f_{coll}$& g & g & $[-20, 20]$\\
$\omega$& 0 & 0 & $[-1, 1]$\\
\hline
\end{tabular}
\label{table:ICBound}
\end{table}
Moreover, if we assume that the final $Y$-position is free, but bounded between [-2, 2], [-5, 5] or [-10, 10] meters, depending on the bounds, the designed trajectories are either Zig-Zag or pure oblique trajectories as illustrated in \reffig{fig:YfreeAllinOne}. For the mentioned bounds, the time of the optimal trajectory is 6.80, 4.60 and 3.39 seconds, respectively. It is obvious that the pure oblique trajectory is faster than the Zig-Zag trajectories, but the latter requires a much smaller space in the 2D space. Also, by decreasing the admissible interval in the horizontal direction, the optimal trajectory time will increase.

\begin{figure}[tbp]
	\centering
	\includegraphics[trim={4.5cm 8cm 4cm 8cm},clip ,width=\linewidth]{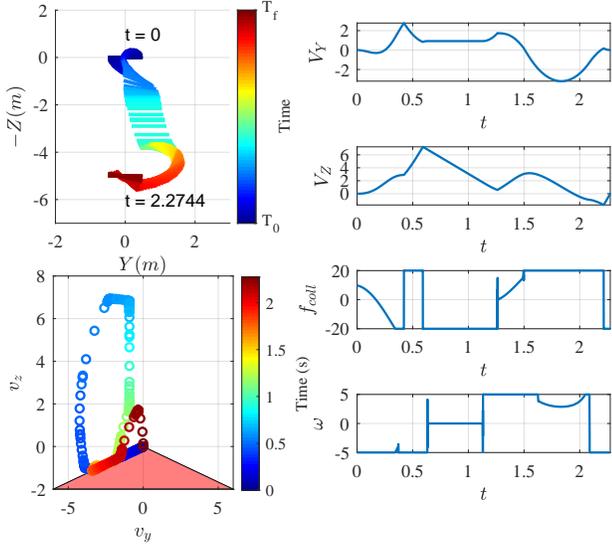}
	\caption{Optimal 2D trajectory design for equal initial and final $Y$, no limitation is considered for $\phi$ allowing for flip maneuvers. It shows that although the final $Y$ is fixed and equal to initial $Y$, the trajectory with flips is faster than ones with simple oblique maneuvers as illustrated in \reffig{fig:YfreeAllinOne}. That is because of the maximum acceleration in the vertical direction is gravitational acceleration, however using flips, the quadcopter can use its thrusts to increase the acceleration in the downward direction.}
	\label{fig:YfixPhifreeAllinOne}
\end{figure}

Also, if we assume that the $\phi$ angle is free, it is permissible that the system has flip maneuvers in its trajectory. With this assumption, the optimal time of the optimal trajectory is 2.27 seconds. \reffig{fig:YfixPhifreeAllinOne} illustrates that the trajectory with flips when the system has freedom in the $\phi$ angle; however it should come back to its initial $Y$ position because of predefined final conditions. Due to the upside down maneuver, the quadcopter can descend with an acceleration higher than the gravity. Moreover, because of this rotation, the descent velocity will be positive in the body fixed frame where there is no velocity constraint in positive body fixed frame $z$ velocity region and the quadcopter can descent freely with high velocity.

\begin{figure}[tbp]
	\centering
	\includegraphics[trim={4.5cm 8cm 4cm 8cm},clip ,width=\linewidth]{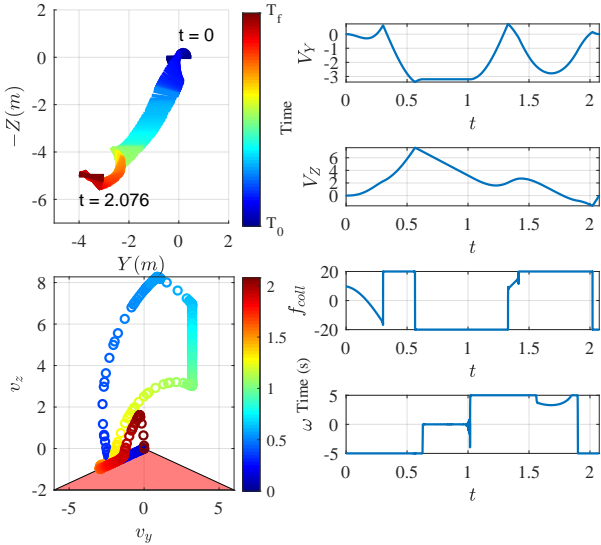}
	\caption{Optimal 2D trajectory design with no limitation on final $Y$, horizontal displacement, and $\phi$ angle interval, allowing for oblique and flips maneuvers, simultaneously. We can see that using oblique maneuver in the combination of flip maneuvers can decrease the final trajectory time. This trajectory is the fastest trajectory that can be designed using the mentioned assumptions, limitations and constraints.}
	\label{fig:YfreePhifreeAllinOne}
\end{figure}

Finally, if we assume that the $Y$ displacement is free as well as the $\phi$ angle, the trajectory would be a combination of the oblique and flip maneuvers, illustrated in \reffig{fig:YfreePhifreeAllinOne}. This trajectory takes 2.07 seconds for descending 5 meters. It shows that the fastest trajectory is a trajectory, including flips and oblique maneuvers. A comparison between methods depending on admissible intervals and allowed maneuvers, is provided in \reftab{table:review}.
\begin{table}[tbp]
	\centering
	\caption{Optimal trajectories with different maneuvers and admissible intervals, to descend 5 meters in the $Z$ direction. Initial and final $\phi$ angles are zero, as well as initial $Y$ condition.}
	\begin{tabular}{ccccccc}
		\hline \noalign{\smallskip} 
		\textbf{Y Intervals} & $\boldsymbol{\phi}$ \textbf{Intervals} & \textbf{Final Y} & \textbf{Ob.} & \textbf{Z.Z.} & \textbf{Fl.} & \textbf{Traj. Dur.} \\
		\hline\noalign{\smallskip}
		 $[-15, 15]$ & $[-\pi/3, \pi/3]$ & fixed(0) & \xmark & \checkmark &  \xmark & 5.33 s\\
		 $[-2, 2]$ & $[-\pi/3, \pi/3]$ & bounded &  \xmark & \checkmark &  \xmark & 6.80 s\\
		 $[-5, 5]$ & $[-\pi/3, \pi/3]$ & bounded &  \xmark & \checkmark &  \xmark & 4.60 s\\
		 $[-10, 10]$ & $[-\pi/3, \pi/3]$ & free &  \checkmark & \xmark &  \xmark & 3.39 s\\
		 $[-15, 15]$ & free & fixed(0) &  \xmark & \xmark &  \checkmark & 2.27 s\\
		 $[-15, 15]$ & free & free &  \checkmark & \xmark &  \checkmark & 2.08 s\\
		\hline
	\end{tabular}
	\label{table:review}
\end{table}

Based on a brief review of different optimal trajectories, provided in \reftab{table:review}, it is evident that the fastest trajectory is a trajectory including flips and oblique maneuvers, when there are no limitations in the horizontal direction. Moreover, the lowest-speed trajectory is the Zig-Zag trajectory when there are strict limitations on the $Y$ displacement. However, if the flip maneuver is not allowed in the trajectory, the simple oblique trajectory will be the fastest solution for a fast descent, considering the velocity constraint.

\subsection{3D Optimal Trajectories}\label{sec:3Danal}
As we have shown in the 2D analysis, the fastest trajectory is the trajectory where $\phi$ and $Y$ are free. However, if aggressive maneuvers are not desired, we should limit the allowed rotation angle of the quadcopter. In the 2D analysis, between the limited $\phi$ trajectories, $Y$ free trajectories are the fastest. However, in real flights, the quadcopter needs to descend fast, mostly we do not have enough space for such large movement in the $Y$ direction, in confined environments, these trajectories are not possible. As we shall see shortly, the trajectories we will obtain in this 3D analysis are helix type trajectories. Note that for helix type trajectories, a centrifugal force is also present. 
%
%

Due to the centrifugal force in helix type trajectory, the quadcopter should have a pitch angle around the $y$ axis of body frame in addition to a roll rotation around x. The thrust vector in the inertial frame is:
\begin{equation}
\textbf{\textit{T}}^{\text{I}}=\textbf{\emph{R}}_\text{b}^\text{I} 
\begin{bmatrix}
0\\ 0\\-T
\end{bmatrix}.
\end{equation}

If the yaw angle is such that the velocity vector belongs to the  xz  plane in body fixed coordinates, then the component of thrust for the centrifugal force can be calculated from:
\begin{equation}
|\textbf{\emph{T}}_\text{cent}| = \left(\sum \textbf{\textit{F}}^\text{I}\right) . \textbf{\emph{e}}_{\text{n}}
\end{equation}
where $\textbf{\emph{T}}_\text{cent}$ is the component of the thrust vector in the centrifugal direction, and $\textbf{\emph{e}}_{\text{n}}$ is the unity vector perpendicular to the trajectory curve, toward to the curvature center. The magnitude of centrifugal force is:
\begin{equation}
\Vert \textbf{\emph{T}}_\text{cent}\Vert= \text{m} V^2/r .
\end{equation}
where $r$ is the radius of the trajectory curve each time. For finding the quadcopter velocity in the body frame, by assuming that the heading of the robot is in the forward direction of motion, we can use the rotational matrices. The wind velocity with respect to the quadcopter body in the inertial frame is:
\begin{equation}
\textbf{\textit{V}}_{\text{W}_\text{rel}}^{I}=\textbf{\textit{W}}-\textbf{\textit{V}},
\end{equation}
where $\textbf{\emph{W}}$ is the wind speed, and $\textbf{\emph{V}}$ is the quadcopter's velocity. To find the wind velocity in the body frame, we can use:
\begin{equation}
\textbf{\textit{V}}_{\text{W}_\text{rel}}^b=\textbf{\emph{R}}_\text{I}^\text{b} \textbf{\textit{V}}_{\text{W}_\text{rel}}^{I},
\end{equation}
where $\textbf{\emph{R}}_\text{I}^\text{b}$ is the rotation matrix from the inertial to the body frame.

\begin{figure}[tbp]
	\centering
	\begin{subfigure}[b]{0.45\textwidth}
		\includegraphics[trim={3.5cm 9cm 4cm 9cm},clip, width=\textwidth]{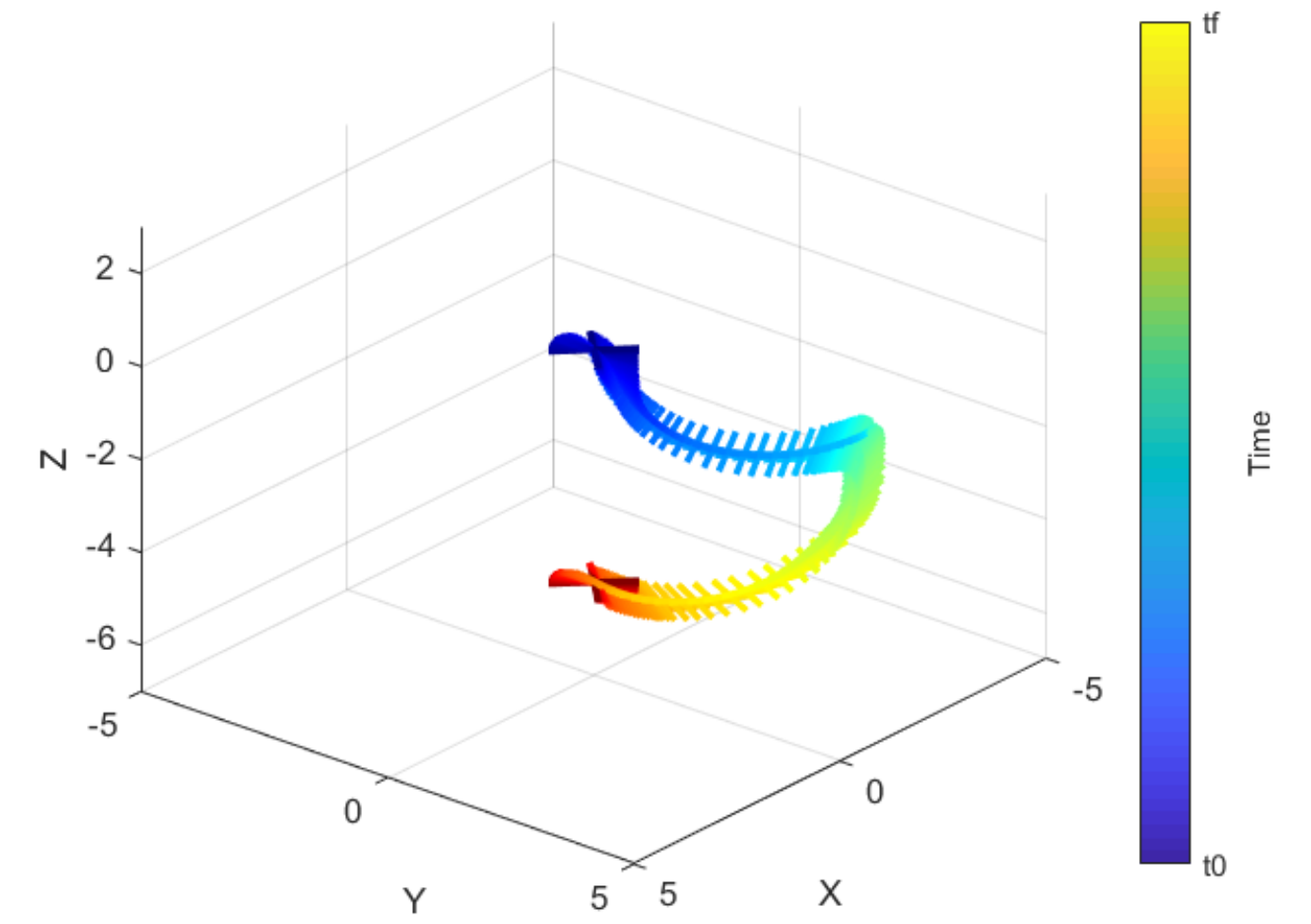}
		\caption{}
		\label{fig:3D-Traj}
	\end{subfigure}
	\begin{subfigure}[b]{0.45\textwidth}
		\includegraphics[trim={4cm 9cm 4.5cm 9cm},clip, width=\textwidth]{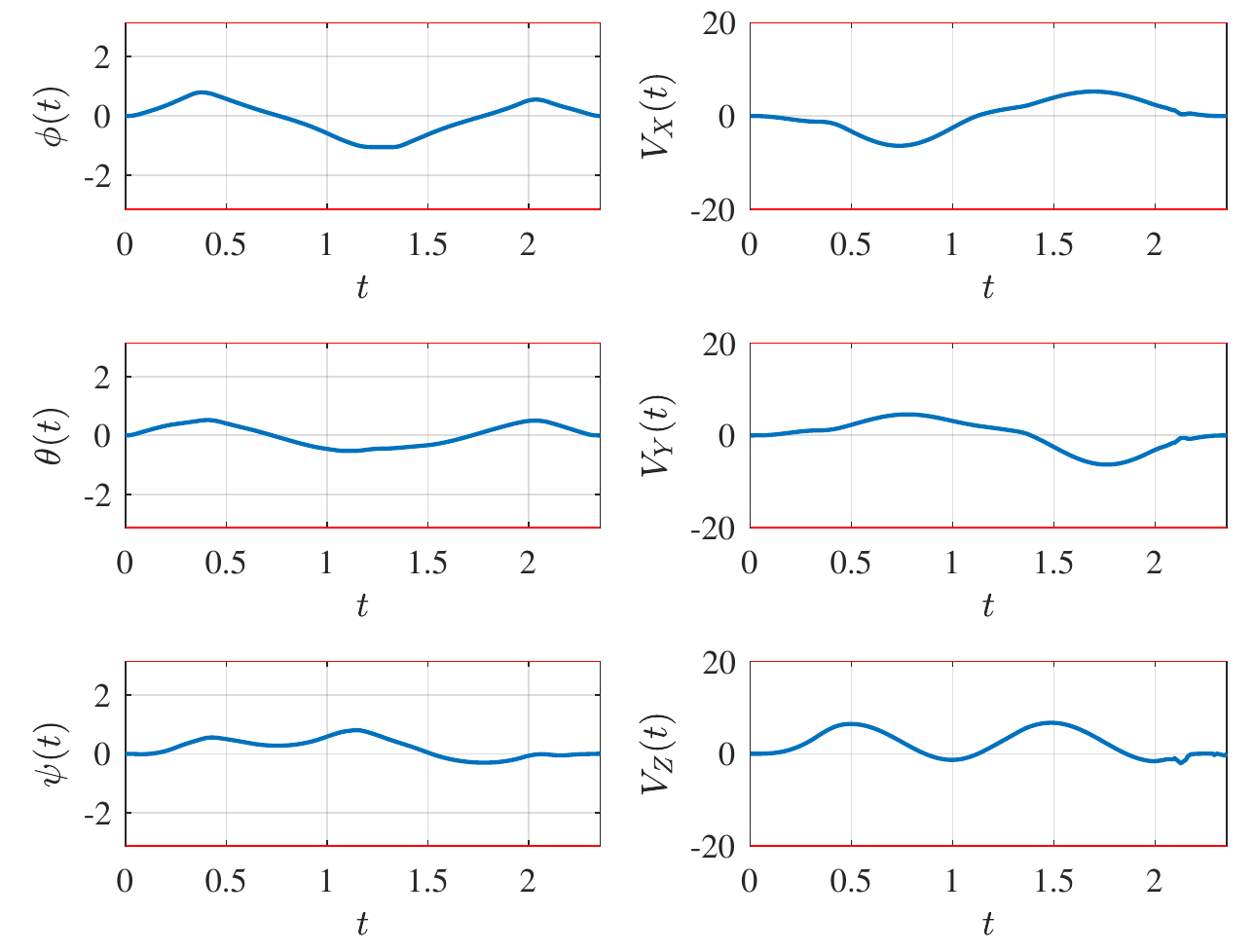}
		\caption{}
		\label{fig:3D-Angles}
	\end{subfigure}
\caption{Optimal 3D trajectory design for fast descent using GPOPS-II, a) Position and 3D Trajectory of generated trajectory, b) Orientation and velocity of the quadcopter in the generated trajectory.
}
\label{fig:3D}
\end{figure}

Using the velocity constraint defined in \refeq{eq:Constraint}, the constraint in the 3D analysis will be:
\begin{equation}
|\textbf{\textit{V}}_{\text{W}_\text{H}}^b|\tan{20^\circ} \geq \textbf{\textit{V}}_{\text{W}_\text{z}}^b,
\end{equation}
where:
\begin{equation}
\left(
\begin{matrix}
V_{\text{W}_\text{x}}^b\\
V_{\text{W}_\text{y}}^b\\
V_{\text{W}_\text{z}}\\
\end{matrix}
\right)
=
\textbf{\emph{R}}_\text{I}^\text{b}
\left(
\begin{matrix}
V_{X}\\
V_{Y}\\
V_{Z}\\
\end{matrix}
\right),
\end{equation}
and:
\begin{equation}
V_{\text{W}_\text{H}}^b=\sqrt{({V_{\text{W}_\text{x}}^b}^2+{V_{\text{W}_\text{y}}^b}^2)}.
\end{equation}
To solve the optimal-time problem for finding the fastest descent trajectory, the cost function in \refeq{eq:cost} is used. The states and control inputs in the 3D analysis are defined as:
\begin{equation}
\left. \begin{array}{ll}
S =(X, V_X, Y, V_Y, Z, V_Z, q, \dot{q}),\\
U = (T_1, T_2, T_3, T_4).\\
\end{array} \right.
\end{equation}
Note that quaternions were used to avoid singularities. Moreover, if we use the equations of motion in 3D space provided in \refsec{subsec:QuadEoM}, considering the \ac{VRS} velocity constraint in \refeq{eq:Constraint}, we will find more complex Hamiltonian conditions that are very difficult to solve analytically. Then, numerical methods are used as for the 2D analysis. Using the GPOPS-II as the numerical solver to find the optimal trajectory, in 3D motion, the helix type trajectory can be generated for minimum time trajectory, illustrated in \reffig{fig:3D}. For this design, 5 meters descent and coming back to its initial $X$ and $Y$ positions are desired. As we have seen in the 2D analysis, for non-flip trajectories, the oblique trajectory is the fastest way to descend. Hence, in 3D analysis, the helix type trajectory is the equivalent of oblique in 2D, with consideration of $Y$-displacement limitation. Then, in real 3D flights, an appropriate non-aggressive (without flips) maneuver for fast descent is a helix type. This trajectory helps motion planners to design a faster descent in their mission, and it is considered for flight experiments.

\section{Experiments}
\label{sec:Experiments}
\begin{figure}[tbp]
	\centering
	\includegraphics[trim={1cm 0.5cm 0cm 0.5cm},clip, width=.9\linewidth]{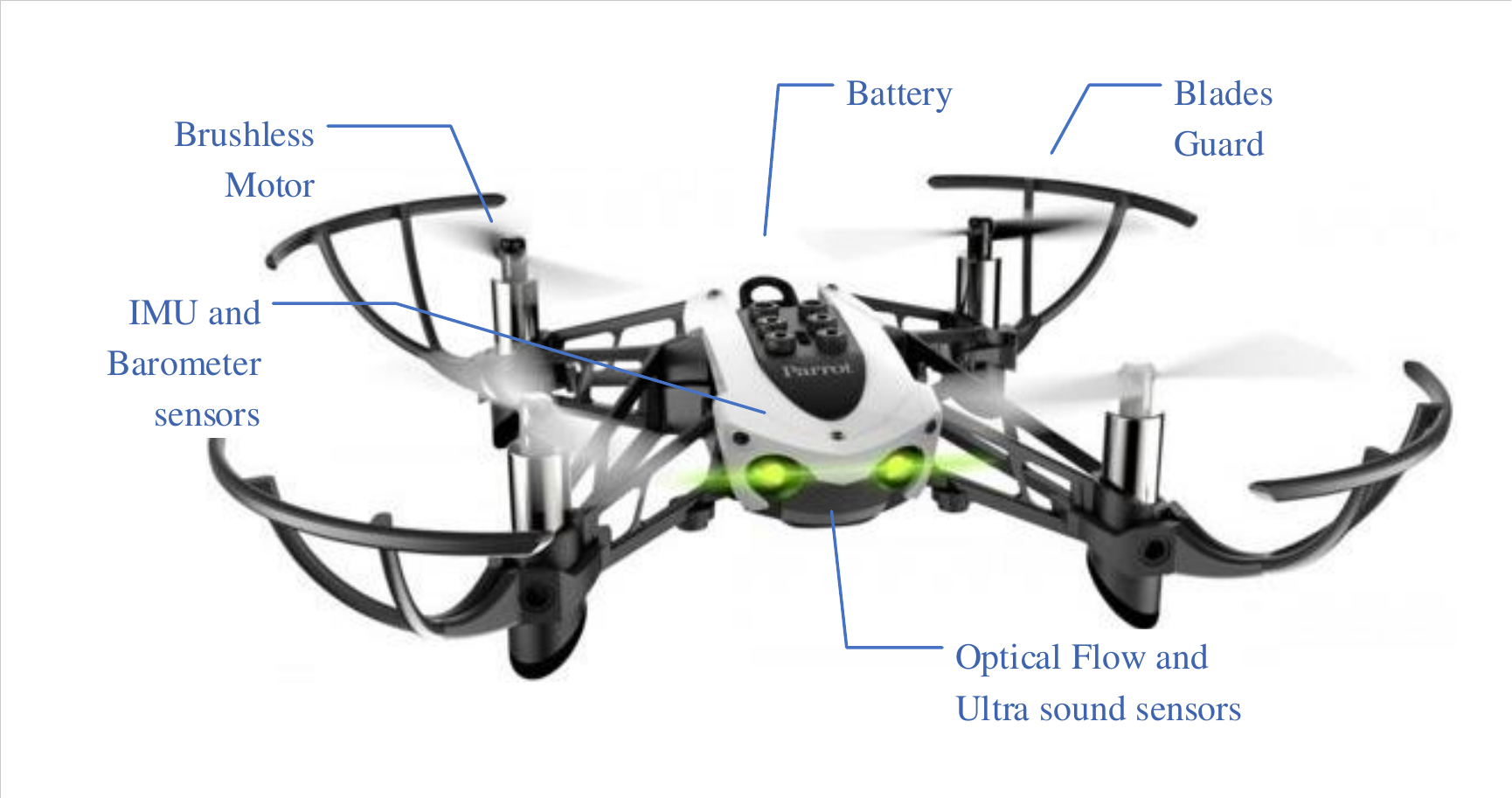}
	\caption{Parrot Minidrone Mambo used for the flight tests}
	\label{fig:Mambo}
\end{figure}

For real flight tests, the Parrot Mambo drone shown in \reffig{fig:Mambo} was utilized. This quadcopter has the following sensors: IMU, barometer, ultrasound, and a downward camera using optical flow algorithm for positioning. The controller and designed trajectories are uploaded on the Mambo from Matlab/Simulink via Bluetooth.

The mass of the Mambo quadcopter is 63 gr, and each blade disk diameter is 6.5 cm. Using the equation~\refeq{eq:InducedVelHover},  the induced velocity at hover is 1.54 m/s.

To show the fluctuations and instability in the \ac{VRS} and \ac{WBS} regions in the real flight, a sinusoidal desired trajectory in the $Z$ direction with increasing frequency is designed, which is given by:
\begin{equation}\label{eq:Zsin}
Z_{des}=A \sin(\omega t) + Z_0,
\end{equation}
where $\omega=\omega_0 t$.

By using reference, the vertical speed would grow slowly, and we can see the effects of entering the prohibited regions, in a limited flight space.  
\begin{figure}[tbp]
	\begin{subfigure}{0.9\linewidth}
		\includegraphics[trim={4cm 9cm 4cm 8.5cm} ,clip, width=\linewidth]{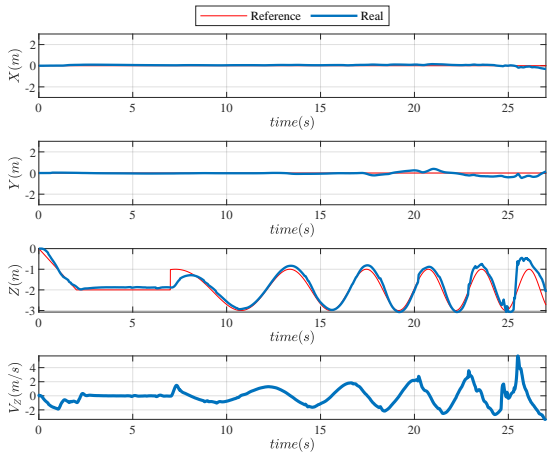}
		\caption{Position and vertical velocity of the flying Quadcopter}
		\label{fig:VarFreq1}
	\end{subfigure}
	\\[7pt]
	\begin{subfigure}{0.9\linewidth}
		\includegraphics[trim={4cm 9cm 4cm 9cm}, clip, width=\linewidth]{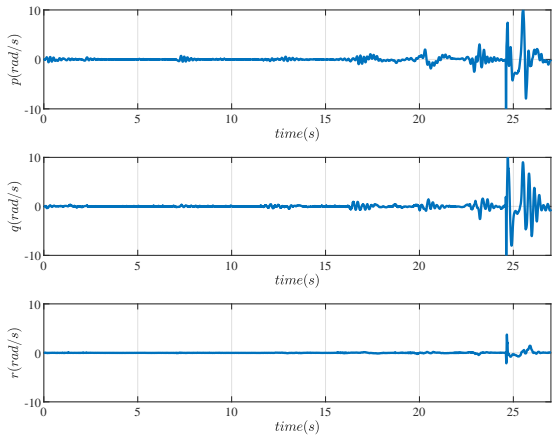}
		\caption{Angular velocity of the quadcopter}
		\label{fig:VarFreq2}
	\end{subfigure}
\caption{Flight tests with variable frequency sinusoidal trajectory in the $Z$ direction that provides increasing descent velocity. This test shows that by entering into the \ac{VRS} or \ac{TWS} regions, the quadcopter experiences fluctuation and by increasing descent velocity, without any horizontal maneuver, the fluctuations are higher till the system becomes completely unstable, at time 24 s.}
\label{fig:VarFreq}
\end{figure}

As we can see in \reffig{fig:VarFreq}, by entering into the \ac{VRS} region (that is 0.6 m/s velocity in the $Z$ direction), there will be low fluctuations in the quadcopter. Entering deeper into the \ac{VRS} and \ac{WBS} regions, more fluctuations are experienced by the quadcopter, till the quadcopter becomes completely unstable, and crashes at 24~seconds. The data in the plots after this moment (24 seconds) pertain to the crashing effects with the ground.

This experiment shows that the \ac{VRS} and \ac{WBS} regions are present in the quadcopter and for the trajectory design for these rotary-wing drones, it is vital to assume the velocity constraint.

To descend as fast as possible, avoiding the \ac{VRS} and \ac{WBS} regions, and assuming non-aggressive maneuvers, a helix type trajectory is suggested from the analysis provided in the previous section. 
\begin{figure}[tbp]
	\begin{subfigure}{0.9\linewidth}
		\centering
		\includegraphics[trim={4.5cm 9cm 4.5cm 8.5cm} ,clip, width=\linewidth]{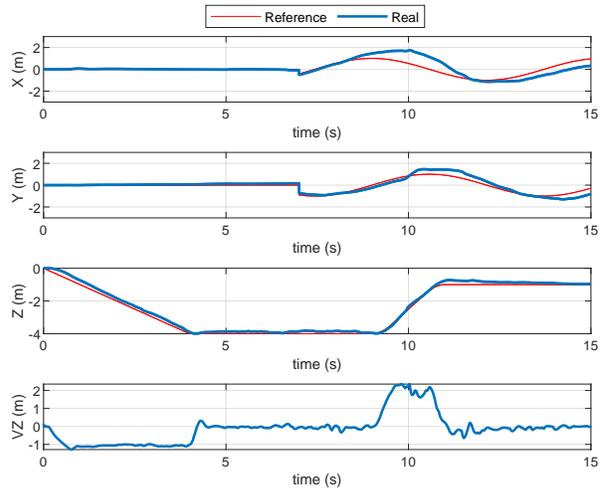}
		\caption{Position and vertical velocity of the flying quadcopter}
		\label{fig:Helix1}
	\end{subfigure}
	\\[7pt]
	\begin{subfigure}{0.9\linewidth}
		\centering
		\includegraphics[trim={4.5cm 9cm 4.5cm 9cm} ,clip, width=\linewidth]{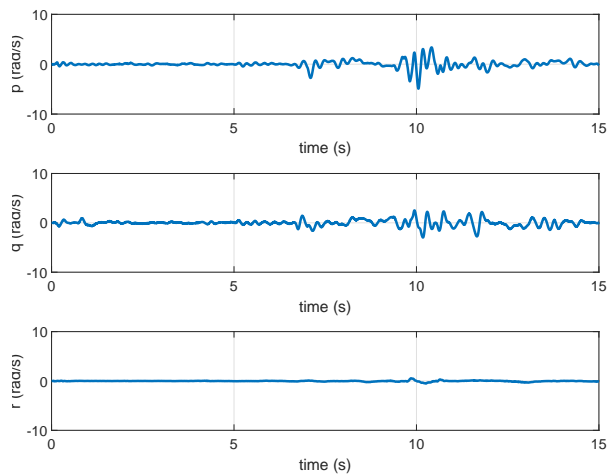}
		\caption{Angular velocity of the quadcopter}
		\label{fig:Helix3}
	\end{subfigure}
	\caption{Drone with a helix type trajectory for achieving the high speeds in descent. It is illustrated that by following the helix type trajectory, the quadcopter can achieve high descent speeds considering the \ac{VRS} and \ac{WBS} constraints on air stream velocities in the blade disk frame.}
	\label{fig:Helix}
\end{figure}
We can see that in \reffig{fig:Helix} even by having descent velocity higher than 0.6 m/s for more than 2 seconds, it can descend fast. Therefore by following a helix type trajectory and thereby adding horizontal speed in the blade disk frame the quadcopter avoids the \ac{VRS} and \ac{WBS} regions and does not become unstable.

In the regular descent, the quadcopter can descent 5 meters with a pure descent in 8.3 seconds. However, using VRS avoiding trajectory, the quadcopter can descend 5 meters in 2.5 seconds which is much less than the pure descent.

\section{Conclusions}
\label{sec:Conclusions}
In this paper, the existence of the \ac{VRS} and \ac{WBS} in the context of quadcopters have been investigated. Subsequently, using experiments in wind tunnel, a normalized model for the quadcopter drones was obtained which is independent from the disk load and the blade disk’s diameters. Then, for trajectory designs, a simple model is provided. Due to the complex optimal problem for minimum time trajectory design, we designed optimal 2D and 3D trajectories for descent, using GPOPS-II package as a numerical solver. Finally, we performed flight experiments which showed that the \ac{VRS} is present for quadcopters. In fact, this might be the cause of instability and crashes in fast descents. Furthermore, it is explained that by increasing the horizontal speed to the blade disk, the fluctuations can be reduced in the flight. By descending with a helix type trajectory as an optimal descending trajectory, the quadcopter can descend fast without entering to the unstable regions.

Possible directions for future work are to use learning methods to identify the unstable boundaries, controlling the quadcopter with learning methods, and investigating the behavior of variable pith quadcopters in the \ac{VRS} and \ac{WBS} regions.

\bibliographystyle{IEEEtran}
\bibliography{Amin_bib}
\vspace{-1cm}
\begin{IEEEbiography}[{\includegraphics[width=\linewidth]{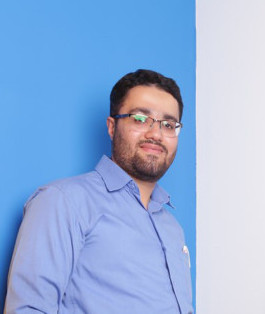}}]{Amin~Talaeizadeh}
	received his B.Sc. and M.Sc. degrees in mechanical engineering from the Amirkabir University of Technology, Tehran, Iran, in 2011 and 2013, respectively. 
	He is currently pursuing the Ph.D. degree from the Sharif University of Technology, Tehran, Iran. He is now a visiting Ph.D. student in the CST group of Eindhoven University of Technology (TU/e). His current research interests include optimal trajectory design, robotics, model-based design, and optimal control.
\end{IEEEbiography}
\vspace{-1cm}
\begin{IEEEbiography}[{\includegraphics[width=\linewidth]{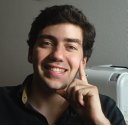}}]{Duarte~Antunes}
received the Licenciatura in Electrical and Computer Engineering in 2005 and his Ph.D. in 2011 from the Instituto Superior Técnico (IST), Lisbon, Portugal,in 2005. From 2011 to 2013 he held a postdoctoral position at the department of mechanical engineering at the Eindhoven University of Technology (TU/e), where he is currently an Assistant Professor.
\end{IEEEbiography}
\vspace{-1cm}
\begin{IEEEbiography}[{\includegraphics[width=\linewidth]{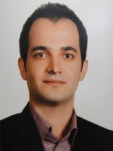}}]{Hossein~Nejat~Pishkenari}
earned his B.Sc., M.Sc. and Ph.D. degrees in Mechanical Engineering from the Sharif University of Technology in 2003, 2005 and 2010, respectively. Then he joined the Department of Mechanical Engineering at the Sharif University of Technology in 2012. Currently, he is directing the Nano-robotics Laboratory and the corresponding ongoing research projects in the multidisciplinary field of Nanotechnology.
\end{IEEEbiography}
\vspace{-1cm}
\begin{IEEEbiography}[{\includegraphics[width=\linewidth]{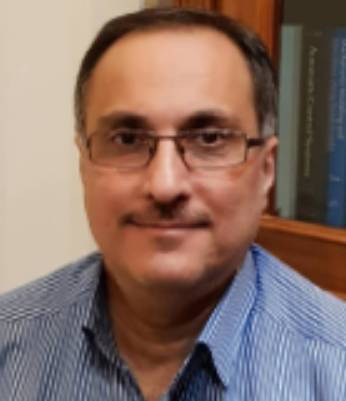}}]{Aria~Alasty}
received his B.Sc. and M.Sc. degrees in Mechanical engineering from Sharif University of Technology (SUT), Tehran, Iran in 1987 and 1989. He also received his Ph.D. degree in Mechanical engineering from Carleton University, Ottawa, Canada, in 1996. At present, he is a professor of mechanical engineering in Sharif University of Technology. He has been a member of Center of Excellence in Design, Robotics, and Automation(CEDRA) since 2001. His fields of research are mainly in Nonlinear and Chaotic systems control, Computational Nano/Micro mechanics and control, special purpose robotics,robotic swarm control, and fuzzy system control.
\end{IEEEbiography}
\end{document}